\documentclass[12pt, a4paper]{article}

\usepackage{amsmath,amssymb}
\usepackage{color}

\numberwithin{equation}{section}
\newcommand{\e}{\epsilon}
\newcommand{\p}{\partial}
\newcommand{\nn}{\nonumber}

\newcommand{\R}{\mathbb{R}}
\newcommand{\Z}{\mathbb{Z}}
\newcommand{\F}{\mathcal{F}}
\newcommand{\V}{\mathcal{V}}
\newcommand{\E}{\mathcal{E}}
\newcommand{\A}{\mathcal{A}}
\newcommand{\T}{\mathcal{S}}
\newcommand{\D}{\mathcal{D}}

\newcommand{\br}{[\ \,,\ ]}
\newcommand{\ow}{\overline{\wedge}}
\newcommand{\Der}{\mathrm{Der}}
\newcommand{\ad}{\mathrm{ad}}

\newcommand{\tF}{\tilde{\F}}
\newcommand{\tV}{\tilde{\V}}
\newcommand{\tT}{\tilde{\T}}
\newcommand{\tE}{\tilde{\E}}
\newcommand{\tp}{\tilde{\p}}
\newcommand{\tu}{\tilde{u}}
\newcommand{\tx}{\tilde{x}}

\newtheorem{thm}{Theorem}[subsection]
\newtheorem{cor}[thm]{Corollary}
\newtheorem{lem}[thm]{Lemma}
\newtheorem{prp}[thm]{Proposition}
\newtheorem{emp}[thm]{Example}
\newtheorem{rmk}[thm]{Remark}
\newtheorem{dfn}[thm]{Definition}

\newenvironment{prf}{\noindent {\it Proof} \ }{\hfill $\Box$}
\newenvironment{prfn}[1]{\noindent {\it Proof of #1} \ }{\hfill $\Box$}

\begin{document}
\title{Jacobi Structures of Evolutionary Partial Differential Equations}
\author{{Si-Qi Liu \quad Youjin Zhang} \\
{\small Department of Mathematical Sciences, Tsinghua University}\\
{\small Beijing 100084, P. R. China} \\
{\small Emails: liusq@mail.tsinghua.edu.cn, youjin@mail.tsinghua.edu.cn}}
\date{\today}
\maketitle

\begin{abstract}
In this paper we introduce the notion of infinite dimensional Jacobi structure to describe the geometrical structure of a class of nonlocal Hamiltonian systems
which appear naturally when applying reciprocal transformations to Hamiltonian evolutionary PDEs. We prove that our class of infinite dimensional Jacobi structures
is invariant under reciprocal transformations. The main technical tool is in a suitable generalization of the classical Schouten-Nijenhuis bracket to the space of
the so called quasi-local multi-vectors, and a simple realization of this structure in the framework of supermanifolds. These constructions are used to the computation
of the Lichnerowicz-Jacobi cohomologies of Jacobi structures. We also introduce the notion of bi-Jacobi structures and consider the integrability of a system of
evolutionary PDEs that possesses a bi-Jacobi structure.
\end{abstract}

\tableofcontents

\section{Introduction}

The notion of Jacobi structure is a
generalization of Poisson structures \cite{lich-1} on finite dimensional smooth manifolds, it
has drawn much attention since its introduction by Lichnerowicz
in 1978 \cite{lich-2}. Jacobi structure unifies several important mathematical structures
which include Poisson manifolds, contact manifolds and locally
conformal symplectic manifolds, it also gives a geometrical description of the local Lie algebra structures
introduced by Kirillov \cite{kil}. On a manifold $M^n$ endowed with a Jacobi structure (which is called a Jacobi manifold)
one can define a bracket on the space of smooth functions $C^\infty(M)$, this bracket generalizes the usual
Poisson bracket defined on a Poisson manifold, although it no longer satisfies the Leibniz
rule, it still satisfies the Jacobi identity, and this is the main property that enable Jacobi structures
to describe certain generalized Hamiltonian structures for dynamical systems.

In the present paper, we study an important class of infinite dimensional
generalized Hamiltonian structures for
evolutionary PDEs, it turns out that infinite dimensional counterpart of Jacobi structures
appear naturally when we try to clarify such class of generalized Hamiltonian structures.

Infinite
dimensional Hamiltonian structures are natural generalizations of
the finite dimensional ones for ordinary differential equations,
they have been playing important role in the study of qualitative
properties of some nonlinear PDEs that arise in fluid mechanics and
other fields of mathematics and physics \cite{arnold}. They are also
one of the main tools in the study of integrability of some
physically important nonlinear evolutionary PDEs, see\cite{faddeev,
gardner,zf, magri, dickey, olver} and references therein. The
most frequently used class of Hamiltonian structures, often with an adjective ``local", for
evolutionary PDEs with unknown functions $u^1,\dots, u^n$ depending
on the spatial variable $x$ and time variable $t$ are given by
Poisson brackets of the form
\begin{equation}\label{ham-zero}
\{F, G\}=\int \frac{\delta F}{\delta u^i(x)} P^{ij} \frac{\delta G}{\delta u^j(x)}\, dx,
\end{equation}
where $F, H$ are local functionals
$$
F=\int f(u;u_x,\dots, u^{(m)}) \, dx,\quad G=\int g(u;u_x,\dots, u^{(k)})\, dx
$$
with densities $f,\ g$ being differential polynomials of the unknown functions, i.e. they are
polynomials in the $x$-derivatives of $u^1,\dots, u^n$ with coefficients depending smoothly on $u^i$, and
$P^{ij}, i,j=1,\dots n$ are differential operators in $\frac{\p}{\p x}$ whose coefficients
are differential polynomials. A Hamiltonian system of evolutionary PDEs associated to a
local functional $H$ have the representation
$$
\frac{\p u^i}{\p t}=\{u^i(x), H\}=P^{ij} \frac{\delta H}{\delta u^j},\quad i=1,\dots, n.
$$
By extending the notion of Schouten-Nijenhuis bracket defined on the
space of multi-vectors of a finite dimensional manifold $M$ to that
of the infinite jet space $J^\infty(M)$
of  $M$, one can represent an
infinite dimensional Hamiltonian structure of the above form in
terms of a local Poisson bi-vector, and can define the associated Poisson
cohomologies \cite{lich-1, olver, DZ}. This leads to an analogue of
the Darboux theorem on the canonical forms of a class of infinite
dimensional Poisson brackets which possess hydrodynamic limits
\cite{Ge, magri-2, DZ}. Such canonical forms play an important role in the
study of the problem of classification of bihamiltonian
integrable hierarchies under Miura type transformations \cite{DZ,
LZ-1, DLZ-1, DLZ-2},  these transformations are certain type of changes of the dependent variables of the hierarchies generalizing the well-known
transformation between the KdV and the modified KdV equations \cite{miura}.

Though the class of local Hamiltonian structures is preserved under the Miura type transformations
considered in \cite{DZ},
the locality of the Hamiltonian structures are in general no longer preserved under changes of the independent variables $x, t$ of the associated evolutionary PDEs,
and in particular under the so called reciprocal transformations of the following form
\begin{equation}
d \tx=\rho\,dx+\sigma\,dt,\quad d \tilde{t}=dt. \label{eq-reci}
\end{equation}
Here $\rho, \sigma$ are differential polynomials of $u^1,\dots, u^n$ which give a conservation law $\frac{\p\rho}{\p t}=\frac{\p\sigma}{\p x}$ of the Hamiltonian system. Nonetheless such reciprocal transformations still preserve the conservation laws
and are important in the study of the system of evolutionary PDEs. A typical
example of their applications to the study of nonlinear integrable systems
is provided by the Camassa-Holm equation \cite{CH, CHH, FF, Fu} and its generalizations\cite{LZ-1, CLZ-1}.
The relation between the Camassa-Holm equation and the first negative flow of the KdV hierarchy
established by a reciprocal transformation plays a crucial role in the study
of integrability of the equation and of the properties of its solutions. Due to the importance of Hamiltonian
structures played in the study of nonlinear evolutionary PDEs, a natural question arises: What is the form of a
nonlocal Hamiltonian structure to which a local Hamiltonian structure is transformed by a reciprocal transformation?

The explanation of the above question will help us to generalize the classification program for bihamiltonian integrable hierarchies under the action of Miura group (consisting of Miura type transformations) to the one
that includes also reciprocal transformations \cite{DZ, LZ-1, DLZ-1}.
The first crucial step toward answering the above question was made in \cite{FP}, where Ferapontov and Pavlov
studied the transformation law of the Hamiltonian structure
of a system of hydrodynamic type
\begin{equation}\label{sys-hyd}
\frac{\p u^i}{\p t}=V^i_j(u) u^j_x,\quad i=1,\dots,n
\end{equation}
under various reciprocal transformations.
One of their results shows that if the system \eqref{sys-hyd} possesses a Hamiltonian structure of hydrodynamic type,
i.e. the Hamiltonian operator $P$ has the form
\begin{equation}
P^{ij}=g^{ij}(u) \p_x+\Gamma^{ij}_k(u) u^k_x,
\end{equation}
then under a reciprocal transformation \eqref{eq-reci},
where $\rho$ and $\sigma$ are smooth function of $u^1, \dots, u^n$, the system
\eqref{sys-hyd} is transformed to a system of hydrodynamic type which possesses a nonlocal Hamiltonian structure given
by an operator of the form
\begin{equation}
\tilde{P}^{ij}=\tilde{g}^{ij}\p_{\tx}+\tilde{\Gamma}^{ij}_ku^k_{\tx}+\tilde{Z}^i_ku^k_{\tx}\p_{\tx}^{-1}u^j_{\tx}+
u^j_{\tx}\p_{\tx}^{-1}\tilde{Z}^j_ku^k_{\tx}. \label{ham-1}
\end{equation}
Ferapontov and Pavlov also considered more general reciprocal transformations which gives more general nonlocal Hamiltonian structures.

To go further in studying the above question along this line, we need to clarify the following points:
\begin{itemize}
\item{} The meaning of actions of the integral operator $\p_{\tx}^{-1}$ on
the variational derivatives of local functionals, and the validity of the Jacobi identity for the Poisson brackets.
\item{} An intrinsic derivation of the transformation rule of a Hamiltonian structure under a reciprocal transformation.
Note that in \cite{FP} the Hamiltonian operators \eqref{ham-1} were derived based on the forms of the transformed system of equations
under reciprocal transformations.
\item{} The transformation rule for more general
Hamiltonian structures with dispersive terms depending on higher order $x$-derivatives of the unknown
functions under reciprocal transformations.
\end{itemize}

The aim of the present paper is to clarify the above points and, in particular, to give a geometrical description of the class of nonlocal
Hamiltonian structures that arise when we apply reciprocal transformations \eqref{eq-reci} to Hamiltonian evolutionary PDEs. Observe
that for a local functional $\tilde{L}$ of $u(\tx)$ we have the
following identity:
\begin{equation}
\p_{\tx}^{-1}\left(u^i_{\tx}\frac{\delta \tilde{L}}{\delta
u^i(\tx)}\right)=-\tilde{E}(\tilde{L}),
\end{equation}
where the differential polynomial $\tilde{E}(\tilde{L})$ is obtained by the action of a linear differential operator
on the density of $\tilde{L}$ and is called the energy of the functional $\tilde{L}$ (see Section 3.2 below).
With this identity we can get rid of the integral operator
that appears in the nonlocal Hamiltonian operator $\tilde{P}$ in \eqref{ham-1} by representing the associated Poisson bracket in
the following form
\begin{align}
\{\tilde{F}, \tilde{G}\}_{\tilde{P}}=&\int \left(\frac{\delta \tilde{F}}{\delta u^i(\tx)}\left(\tilde{g}^{ij}\p_{\tx}+\tilde{\Gamma}^{ij}_ku^k_{\tx}\right)
\frac{\delta \tilde{G}}{\delta u^j(\tx)}\right.\nn\\
&\left.\quad+\tilde{Z}^i_ku^k_{\tx}\left(\tilde{E}(\tilde{F})\frac{\delta \tilde{G}}{\delta u^i(\tx)}-\frac{\delta \tilde{F}}{\delta u^i(\tx)}\tilde{E}(\tilde{G})
\right)\right)d\tx.\label{ham-e}
\end{align}
In this expression the operator $\p_{\tx}^{-1}$ does not appear, hence the above bracket is in fact a local object.

The bracket \eqref{ham-e} is an alternating bilinear map from the space of local functionals to itself, and generalizes the Poisson bracket \eqref{ham-zero}
of the usual form. Our approach of studying the nonlocal Hamiltonian structures is based on properties of the space
of alternating multi-linear maps
(we call them generalized multi-vectors) from the space of local functionals to itself. On this space there is defined a natural bilinear operation called the
Nijenhuis-Richardson bracket \cite{NR-1, NR-2}. Restricting it to an appropriate subspace consisting of the so called quasi-local
multi-vectors we obtain a Schouten-Nijenhuis
bracket which enables one to represent a nonlocal Hamiltonian structure mentioned above
by a quasi-local bivector whose Schouten-Nijenhuis bracket with itself vanishes.

A quasi-local bivector $J$ can be represented by a local bivector $\Lambda$ and a local vector $X$.
In the case when the degrees of $\Lambda$ and $X$ are equal to zero, $\Lambda$ becomes a bivector field  and $X$ a vector field on the manifold $M$,
and the vanishing of the Schouten-Nijenhuis bracket of $J$ with itself is equivalent to the following conditions
\[
[\Lambda,\Lambda]=2 X\wedge\Lambda,\quad [X, \Lambda]=0,
\]
here $\br$ is the usual Schouten-Nijenhuis bracket on $M$, so $J$ is in fact a Jacobi structure on $M$.
Thus quasi-local Hamiltonian structures are just infinite dimensional generalizations of Jacobi structures on the finite dimensional manifold,
so we also call them Jacobi structures. It is interesting to note  that the reciprocal transformations of Jacobi structures of degree zero correspond to the conformal
changes of Jacobi structures in the finite dimensional case. 

The main content and the organization of the paper is as follows.

In Section 2, we first recall the definition of the algebra of differential polynomials, introduce some operators on this algebra that will be used in subsequent
sections and prove some important identities that are satisfied by these operators. We then define the space of quasi-local multi-vectors on the infinite jet space
$J^{\infty}(M)$ and prove the existence of a Schouten-Nijenhuis bracket on the space of quasi-local multi-vectors. The main theorem of this section is Theorem
\ref{thm-main} which gives a simple expression of the Schouten-Nijenhuis bracket in terms of super variables. We also give the transformation rule of the quasi-local
multi-vectors under Miura type transformations and reciprocal transformations, the main results are contained in Theorem \ref{thm-miura}, \ref{thm-miura2}
and Theorem \ref{thm-reci}, \ref{thm-reci2}. 

In Section 3, we first give the definition of infinite dimensional Jacobi structures for systems evolutionary PDEs, and show that the degree zero and degree one Jacobi
structures coincide respectively with the finite dimensional Jacobi structure on $M$ and the non-local Hamiltonian structures of hydrodynamic type that were studied
in \cite{FP}. One important property of the degree one Jacobi structures, call the Jacobi structures of hydrodynamic type, is that they can be converted to local
Hamiltonian structures under certain reciprocal transformations, see Theorem \ref{zh-jht}. Along this line, we proceed to give a necessary and sufficient condition
under which a Jacobi structure is transformed to a local Hamiltonian structure by a given reciprocal transformation. We then prove an analogue of Darboux theorem
for Jacobi structures by proving the vanishing of some Lichnerowicz-Jacobi cohomologies, which states that any deformation of a Jacobi structure of hydrodynamic
type is trivial under Miura type transformations and reciprocal transformations. We also consider properties of a system of evolutionary PDEs that possesses a
{\it bi-Jacobi structure}, which is an analogue of a bihamiltonian system. By using the triviality of Lichnerowicz-Jacobi cohomologies, we show the existence,
for generic cases, of an infinite number of commuting flows and conservation laws for the system, and thus give evidence of integrability of such a system.
It is also proved that, for a given bihamiltonian structure, if a reciprocal transformation keeps its locality, then its central invariants (which are defined
in \cite{LZ-1, DLZ-1}) are also preserved after the transformation, see Theorem \ref{thm-civ0}.

In Section 4, we give some examples of systems of evolutionary PDEs that possess Jacobi structures. The first example is the KdV equation and the Camassa-Holm equation
each of which possesses (at least) four Jacobi structures. We show that a reciprocal transformation which transforms the Camassa-Holm hierarchy to the negative
flows of the KdV hierarchy establishes the relation of the associated Jacobi structures. The second example is about a special reciprocal transformation that
acts on the bihamiltonian structure defined on the jet space of a Frobenius manifold, we shown that such reciprocal transformations preserve the locality of the
bihamiltonian structures, and that they correspond to the inversion symmetries of the WDVV equations that are given in \cite{D-1} by Dubrovin.   

We end the paper with some concluding remarks in Section 5.

\section{Quasi-local multi-vectors}

\subsection{The differential polynomial algebra $\A$}\label{sec-2}

Let $M$ be a smooth manifold of dimension $n$, $J^\infty(M)$ be the infinite jet space of $M$. Recall that $J^\infty(M)$ is a fiber bundle with fiber $\R^\infty$
which is the projective limit of the projective system
\[\left(\left\{\R^{kn}\right\}_{k\ge1},\ \left\{\pi_{k,l}:\R^{kn} \to \R^{ln}\right\}_{k \ge l}\right),\]
where the projection $\pi_{k,l}$ is the forgetful map
\[\pi_{k,l}: (u^{(1)}, u^{(2)}, \dots, u^{(l)}, \dots, u^{(k)}) \mapsto (u^{(1)}, u^{(2)}, \dots, u^{(l)})\]
with $u^{(k)}=(u^{1,k},\dots,u^{n,k})\in\R^n$. Let $U\times \R^\infty,\, V\times \R^\infty$ be two charts of a local trivializations of $J^\infty(M)$ with local coordinates
\[U:\, (u^1,\dots,u^n;\,u^{(1)},\,u^{(2)},\,\dots),\quad V:\, (v^1,\dots,v^n;\, v^{(1)},\, v^{(2)},\,\dots),\]
then the transition functions of $J^\infty(M)$ are given by the chain rule of higher order derivatives:
\begin{equation}
v^{i,1}=\frac{\p v^i}{\p u^j}u^{j,1},\
v^{i,2}=\frac{\p v^i}{\p u^j}u^{j,2}+\frac{\p^2 v^i}{\p u^k \p u^l}u^{k,1}u^{l,1},\ \dots. \label{tran}
\end{equation}
We remark that the bundle $J^\infty(M)$ is not a vector bundle, since the structure group is not the general linear group of the fiber.

A function $f\in C^\infty(J^\infty(M))$ is called a differential polynomial if it depends on the jet variables polynomially in certain local coordinate system.
All differential polynomials form a subalgebra of $C^\infty(J^\infty(M))$, we denote this subalgebra by $\bar{\A}$.
Note that this definition does not depend on the choice of local coordinate system due to the form of the transition functions \eqref{tran}.
On the ring $\bar{\A}$ there is a gradation defined by
\begin{equation}
\deg u^{i,s}=s,\ \deg f(u)=0\ \mbox{ for } f\in C^\infty(M).
\label{grad}
\end{equation}
We denote by $\A$ the completion of $\bar{\A}$ w.r.t. this gradation, and call it the differential polynomial ring of $M$.
By abusing notations we also call elements of $\A$ differential polynomials, though they
may be infinite sums of differential polynomials (as defined above) with increasing degrees.

There is a vector field globally defined on $J^\infty(M)$ which is defined as 
\begin{equation}
\p=\sum_{i=1}^n\sum_{s=0}^\infty u^{i,s+1}\frac{\p}{\p u^{i,s}}, \label{deri}
\end{equation}
where $u^{i,0}=u^i$. This vector field yields a derivation on $\A$ in an obvious way.

We assume henceforth that $M$ is a contractible manifold with a local chart $U:\, (u^1, \dots, u^n)$, and all the computations are done on this chart.
We denote $\p_{i,s}=\frac{\p}{\p u^{i,s}}$ if $s\ge0$, and $\p_{i,s}=0$ if $s<0$.

In the remaining part of this section, we introduce several useful differential operators on $\A$, and prove some important properties of them.
More details on the differential operators on $\A$ can be found in the next subsection.

\begin{dfn}
For $i=1, \dots, n$ and integers $\alpha, s\ge0$, the $i$-th higher generalized momentum operator of type $(\alpha,s)$ is defined as
\begin{equation}
p_{i,\alpha,s}=\sum_{t\ge0}(-1)^t \binom{t+s}{s} \p^t\p_{i,\alpha+s+t}\ : \A \to \A.
\end{equation}
When $s=-1$, we denote $p_{i,\alpha, -1}=\p_{i,\alpha-1}$.
\end{dfn}

When $s=0$, the operator $p_{i,\alpha,0}$ is just the generalized momentum operator in the Lagrangian mechanics on jet bundles.
When $\alpha=0$, the operator $p_{i,0,s}$ coincides with the higher Euler operator $\delta_{i,s}$ introduced by Kruskal, Miura, Gardner, and Zabusky
\cite{KMGZ} (see also \cite{Ge}). In particular, if $\alpha=s=0$, the operator $p_{i,0,0}$ is just the $i$-th variational derivative $\delta_i$.

\begin{dfn}
For any integer $s\ge-1$, the higher energy operator is defined as
\begin{equation}
E_s=\sum_{\alpha\ge1} u^{i,\alpha}p_{i,\alpha,s}.
\end{equation}
In particular, $E_{-1}=\p$. The energy operator is defined as $E=E_0-1$.
\end{dfn}

In Lagrangian mechanics, if we regard $L\in\A$ as a Lagrangian,
then $E(L)$ is just the Hamiltonian of $L$, so we call $E$ the {\it energy operator}.

\begin{lem}\label{lem-idt1}
The operators $p_{i,\alpha,s}$, $E_s$, $E$ satisfy the following identities:
\begin{align*}
\mbox{i) }   & p_{i,\alpha,s}\,\p=p_{i,\alpha,s-1},\ E_s\,\p=E_{s-1},\ E\,\p=0,\\
\mbox{ii) }  & \p\,p_{i,\alpha,s}=p_{i,\alpha,s-1}-p_{i,\alpha-1,s},\ \p\,E_s=E_{s-1}-u^{i,1}\delta_{i,s},\ \p\,E=-u^{i,1}\delta_i,\\
\mbox{iii) } & p_{i,\alpha,s}(f \cdot g)=\sum_{t\ge0}(-1)^t\binom{t+s}{s}\left(p_{i,\alpha,s+t}(f)\cdot\p^t(g)+\p^t(f)\cdot p_{i,\alpha,s+t}(g)\right),\\
       & E_s(f \cdot g)=\sum_{t\ge0}(-1)^t\binom{t+s}{s}\left(E_{s+t}(f)\cdot\p^t(g)+\p^t(f)\cdot E_{s+t}(g)\right),\\
       & E(f \cdot g)=\sum_{t\ge0}(-1)^t\left(E_t(f)\cdot\p^t(g)+\p^t(f)\cdot E_t(g)\right)-f \cdot g,\\
\mbox{iv) }  & \sum_{t\ge0}\binom{t+s}{s}\p^t p_{i,\alpha,s+t}=\p_{i,\alpha+s},\\
       & \sum_{t\ge0} \p^t\left(f \cdot p_{i,\alpha,t}(g)\right)=\sum_{t\ge0}\p^t(f)\cdot \p_{i,\alpha+t}(g),\\
       & \sum_{s\ge0}\binom{s+t}{t}\p^sE_{s+t}=\sum_{\alpha\ge1}\binom{\alpha+t}{t+1}u^{i,\alpha}\p_{i,\alpha+t},\\
       & \sum_{t\ge0}\p^t\left(f\cdot E_t(g)\right)=\sum_{t\ge0}\sum_{\alpha\ge1}\p^t\left(f\,u^{i,\alpha}\right)\p_{i,\alpha+t}(g).
\end{align*}
\end{lem}

\begin{prf}
All these identities can be proved by direct computations. We omit the details here.
\end{prf}

Besides these identities, the operators $p_{i,\alpha,s}$, $E_s$, $E$ also satisfy many other identities which are not easy to prove directly.
The following lemma is quite useful when dealing with such identities.

\begin{lem}\label{lem-3}
Let $c:\A\to\A$ be a differential operator, if $c(\A)\subset\R$ then $c=0$.
In particular, if $\p\cdot c=0$, then $c=0$. Here we view $\R$ as the subalgebra of $\A$ consisting of constant functions.
\end{lem}
The proof of the lemma is a little bit technical, we give it in Section \ref{app-a}.

\begin{lem}\label{lem-idt3}
The operators $\delta_{i,t}$, $\delta_j$, $E_t$, $E$ satisfy the following identities:
\begin{align}
\mbox{i) }   & \delta_{i,t}\,\delta_j=(-1)^t\p_{j,t}\,\delta_i,\\
\mbox{ii) }  & E\,\delta_j=\p_{j,0}\,E,\quad E_t\,\delta_j=(-1)^t \p_{j,t}\,E\ (t\ge1),\\
\mbox{iii) } & E_t\,E=(-1)^t\sum_{s\ge0}\binom{s+t}{t}\p^s\,E_{s+t}\,E.
\end{align}
\end{lem}
\begin{prf}
The identity i) is just the closeness condition for a variational one-form.
It can be proved by identifying the two sides with the following operator
\[(-1)^t\sum_{p,q\ge0}(-1)^p\binom{p}{t-q}\p^{p+q-t}\p_{i,p}\p_{j,q}.\]

For the identity ii), note that
\[\p\,E\,\delta_j=-u^{i,1}\delta_i\delta_j=-u^{i,1}\p_{j,0}\delta_i=\p\,\p_{j,0}E,\]
so we have $E\,\delta_j=\p_{j,0}E$ according to Lemma \ref{lem-3}. Next, for the $t=1$ case,
\[\p\,E_1\,\delta_j=\left(E_0-u^{i,1}\delta_{i,1}\right)\delta_j=\p_{j,0}E+\delta_j+u^{i,1}\p_{j,1}\delta_i=\p\left(-\p_{j,1}E\right).\]
The $t\ge2$ cases can be proved similarly by induction on $t$.

The identity iii) can be proved as follow.
\begin{align*}
&(-1)^t E_tE=(-1)^tE_{t+1}\p\,E=(-1)^{t+1}E_{t+1}\left(u^{i,1}\delta_i\right)\\
=&(-1)^{t+1}\sum_{p\ge0}(-1)^p \binom{p+t+1}{t+1}\left(E_{p+t+1}\left(u^{i,1}\right)\p^p\delta_i+u^{i,p+1}E_{p+t+1}\delta_i\right)\\
=&(-1)^{t+1}\sum_{p\ge0}(-1)^p \binom{p+t+1}{t+1}u^{i,p+1}(-1)^{p+t+1}\p_{i,p+t+1}E\\
=&\sum_{\alpha\ge1}\binom{\alpha+t}{t+1}u^{i,\alpha}\,\p_{i,\alpha+t}E
=\sum_{s\ge0}\binom{s+t}{t}\p^sE_{s+t}E,
\end{align*}
here we used $E_{p+t+1}\left(u^{i,1}\right)=0$ and some identities given in Lemma \ref{lem-idt1}.
\end{prf}

\vskip 1em

The following lemma is well-known in the theory of variational calculus, we give an alternative proof here in terms of notations we used above.

\begin{lem}\label{lem-1}
Let $f$ be a nonzero element of $\A$, then $f\in\p\A$ if and only if $f\notin\R$ and $\delta_if=0$ for all $i=1, \dots, n$.
\end{lem}
\begin{prf}
Since $\delta_i\p=0$, the necessary condition is obvious, we only need to prove the sufficient condition.
Without loss of generality, we assume that $f$ is a homogeneous element of $\A$.
so we can set $f=f(u, u^{(1)}, \dots, u^{(N)})\notin\R$.
Introduce the notation $Z_{(i,s)}=\p_{i,s}Z$. Then we have
\[(\delta_i f)_{(j,2N)}=(-1)^N f_{(i,N)(j,N)}=0,\]
so $f$ is linear in all $u^{i,N}$ and we can assume that
\[f=\sum_{i=1}^n f_i(u, u^{(1)}, \dots, u^{(N-1)})\,u^{i,N}+f_0(u, u^{(1)}, \dots, u^{(N-1)}).\]
{}From the identity
\[(\delta_i f)_{(j,2N-1)}=(-1)^N \left(f_{(i,N)(j,N-1)}-f_{(i,N-1)(j,N)}\right)=0\]
it follows the existence of $g(u, u^{(1)}, \dots, u^{(N-1)})\in\A$ such that
\[f_i=g_{(i,N-1)},\]
thus $f-\p\,g$ is a differential polynomial that does not depend on $u^{i,N}$, then the proof of the lemma is
finished by
induction on $N$.
\end{prf}

We introduce the notation $\sim$ as follows:
\[A \sim B,\ \mbox{if } A-B\in\p\A,\]
where $A, B\in\A$. The following lemma is frequently used in this paper.

\begin{lem}\label{lem-2}
Let $X^i\in\A\ (i=1, \dots, n)$, if for any $f\in\A$, we have
\[X^i\delta_i(f)\sim0,\]
then there exists $c\in\R$ such that $X^i=c\,u^{i,1}$.
\end{lem}
\begin{prf}
Denote $Z_j=\delta_j\left(X^i\,\delta_i f\right)$, then $Z_j=0$.

Without loss of generality, we assume that $X^i\ (i=1, \dots, n)$ are homogeneous elements of $\A$.
Let $X^i=X^i(u, u^{(1)}, \dots, u^{(N)})$ which depends nontrivially on $u^{(N)}$. By taking $f=u^i$ we see that $X^i\sim0$,
so $X^i$ depends on $u^{(N)}$ linearly and can be put into the form
\[X^i=\sum_{j=1}^n X^i_j(u, u^{(1)}, \dots, u^{(N-1)})u^{j,N}+X^i_0(u, u^{(1)}, \dots, u^{(N-1)}).\]

If $N=2p$ is even, then by putting $f=(-1)^p\left(u^{i,p}\right)^2/2$ we get
\[Z_{i,(k,2N)}=X^i_k+\delta^{ik}X^i_i=0,\]
so $X^i_j=0$, this contradicts our assumption on $X^i$.

Now let $N$ be odd. We first assume that $N=2p+1\ge3$. By taking  $f=\frac{(-1)^{p+1}}2\left(u^{i,p+1}\right)^2$
we obtain
\[Z_{i,(k,2N+1)}=X^i_k-\delta^{ik}X^i_i=0,\]
it follows that $X^i_k=0$ if $k \ne i$. Next, from
\[Z_{i,(i,2N),(k,N)}=\left(1+2\delta^{ik}\right)X^i_{i,(k,N-1)}=0\]
we know that $X^i_i$ does not depend on $u^{(N-1)}$. Then we have
\[Z_{i,(i,2N)}=2X^i_{0,(i,N-1)}+\p X^i_i=0.\]

Now let us take $f=\frac{(-1)^{p+1}}6 \left(u^{i,p+1}\right)^3$, then from
\[Z_{i,(i,2N)}=u^{i,p+1}\left(2X^i_{0,(i,N-1)}+\p X^i_i\right)-N\,X^i_i\,u^{i,p+2}=0\]
it follows that $X^i_i=0$. This contradicts our assumption that $X^i$ depends nontrivially on $u^{(N)}$,
so we must have $N\le1$.

If $N=1$, similar to the case when $N$ is old and $N\ge 3$, we can obtain
\[X^i=c^i\,u^{i,1},\]
where $c^i\in\R$. Finally we take $f=f(u)$, then
\[Z_j=\sum_{i=1}^n (c^i-c^j)f_{(i,0)(j,0)} u^{i,1}=0 ,\]
so $c^i=c^j$. The lemma is proved.

Finally, we remark that when $X^i=c\,u^{i,1}$,
\[X^i\,\delta_i(f)=c\,u^{i,1}\delta_i(f)=\p\left(-c\,E(f)\right)\sim0,\]
so the condition of the lemma is indeed satisfied.
\end{prf}

The following two lemmas are very useful when computing $\delta_i$ and $E$.

\begin{lem}\label{lem-6}
Let $D$ be a first order differential operator over $\A$, if \[D\p=\p\left(D-D(1)\right),\]
then for any $f\in\A$, we have
\[D(f)\sim X^i\delta_i(f)-a\,E(f),\]
where $a=D(1)$, $X^i=\left(D-D(1)\right)(u^i)$.
\end{lem}
\begin{prf}
The general form of $D$ is
\[D=a+\sum_{i,s}X^{i,s}\p_{i,s},\]
where $X^{i,0}=X^i$. The condition $D\p=\p\left(D-D(1)\right)$ is equivalent to
\[X^{i,s+1}=\p\left(X^{i,s}\right)-a\,u^{i,s+1},\]
which implies
\[X^{i,s}=\p^s \left(X^i\right)-\sum_{t=1}^s\p^{s-t}\left(a\,u^{i,t}\right).\]
Then by using the definition of $\delta_i$ and $E$, the lemma is proved immediately.
\end{prf}

\begin{lem}\label{lem-7}
Let $f\in\A$ be a differential polynomial, if for any $g\in\A$, $f\,g\sim0$, then we have $f=0$.
\end{lem}
\begin{prf}
First we take $g=1$, then there exists $h\in\A$ such that $f=\p(h)$. Next we take $g=f$, so we have $\delta_i(f^2)=0$.
Without loss of generality, we assume that $h$ is a homogeneous differential polynomial, so there exists $N\in\mathbb{N}$ such that
$h_{(i,s)}=0$ for any $s>N$.
On the other hand, we have
\[0=\left(\delta_i(f^2)\right)_{(i,2N+2)}=(-1)^{N+1}\left(f^2\right)_{(i,N+1)(i,N+1)}=2(-1)^{N+1}\left(h_{(i,N)}\right)^2,\]
so $h$ must be a constant, hence $f=0$. The lemma is proved.
\end{prf}

\begin{cor}\label{lem-8}
Let $f\in\A$ be a differential polynomial, and $g_1, \dots, g_n, h\in\A$.
If for any first order differential operator $D$ satisfying $D\p=\p\left(D-D(1)\right)$,
\[D(f)\sim X^i\,g_i-a\,h,\]
where $a=D(1)$, $X^i=\left(D-D(1)\right)(u^i)$, then we have $g_i=\delta_i(f)$ and $h=E(f)$.
\end{cor}

\subsection{Differential operators on $\A$} \label{app-a}

In this subsection, we define the algebra of differential operators over $\A$, and prove Lemma \ref{lem-3}.
The readers who are not interested in the proof of Lemma \ref{lem-3} can skip this subsection.

\begin{dfn}
We define the following algebras
\begin{align*}
&\D_M=C^\infty(M)\otimes\R[\p_{i,0}|i=1, \dots, n],\\
&\D_M'=\A\otimes\R[\p_{i,0}|i=1, \dots, n],\\
&\bar{\D}=\A\otimes\R[\p_{i,s}|i=1, \dots, n,\ s\ge0],
\end{align*}
whose elements are called differential operators over $M$, differential operators over $M$ with coefficients in $\A$,
and {\em bounded} differential operators over $\A$ respectively. 

The algebra $\bar{\D}$ possesses a gradation
\[\bar{\D}=\bigoplus_{d\in\Z}\bar{\D}_d,\ \bar{\D}_d=\{\sum_{I,S}f^{I,S}\p_{I,S}\,\left|\,\deg f^{I,S}-|S|=d\right.\},\]
which induces a decreasing filtration
\[\bar{\D}\supset\cdots\supset\bar{\D}_{(d-1)}\supset\bar{\D}_{(d)}\supset\bar{\D}_{(d+1)}\supset\cdots,\ \bar{\D}_{(d)}=\sum_{i \ge d}\bar{\D}_i,\]
we denote the topological completion of $\bar{\D}$ w.r.t. the above filtration by $\D$, and call it the differential operator algebra over $\A$.
\end{dfn}
Here we used the multiple index notations:
\[(I,S)=(i_1, s_1; \cdots; i_m, s_m),\ \p_{I,S}=\p_{i_1, s_1}\cdots\p_{i_m, s_m},\]
where $m=0, 1, 2, \dots$, $|I|=m$, and $|S|=s_1+\cdots+s_m$.

Before proving Lemma \ref{lem-3}, we need some preparations.

\begin{lem}\label{alem-1}
Let $D\in\D_M$, if for any $f\in C^\infty(M)$, $D(f)\in\R$, then $D=0$.
\end{lem}
\begin{prf}
The general form of $D$ reads
\[D=\sum_{|K|<m} a_K\p^K,\]
where $K=(k_1, \dots, k_n)$ is the multiple index, and
\[\p^K=\p_{1,0}^{k_1}\cdots\p_{n,0}^{k_n},\ |K|=k_1+\cdots+k_n.\]
By choosing sufficiently many linearly independent $f_1, \dots, f_N$ and by using the property of Wronskian determinants,
one can prove that all the coefficients $a_K$ vanish. The lemma is proved.
\end{prf}

\begin{lem}\label{alem-2}
Let $D\in\D_M'$, if for any $f\in C^\infty(M)$, $D(f)\in\R$, then $D=0$.
\end{lem}
\begin{prf}
Suppose
\[D=\sum_{I,S}u^{I,S}D_{I,S},\]
where $(I,S)=(i_1, s_1; \cdots; i_m, s_m)$ with $s_i\ge1$ and$D_{I,S}\in\D_M$, then every $D_{I,S}$ satisfies the condition of Lemma \ref{alem-1}, so $D_{I,S}=0$, then $D=0$. The lemma is proved.
\end{prf}

\vskip 1em

\begin{prfn}{Lemma \ref{lem-3}}
According to the definition of $\D$, every $c\in\D$ can be decomposed into
\[c=\sum_{m\ge0}c_{m},\ c_{m}=\sum_{|I|=m}\sum_{S}c_{m}^{I,S}\p_{I,S},\]
where $(I,S)=(i_1, s_1; \cdots; i_m, s_m)$ with $s_i\ge1$, $c_{m}^{I,S}\in \D_M'$.

First we take $f\in C^\infty(M)$, so
\[c(f)=c_{0}(f),\]
since $c(f)\in\R$, according to Lemma \ref{alem-2}, we have $c_{0}=0$.

Next, let $m$ be a positive integer, assume that for any multiple index $(I',S')$ with $|I'|=m'<m$ we have proved $c_{m'}^{I',S'}=0$,
then for a multiple index $(I,S)$ with $|I|=m$, let $f=\tilde{f}\,u^{I,S}$ where $f\in C^\infty(M)$, then we have
\[c(f)=C(I,S)\,c_{m}^{I,S}(\tilde{f}),\]
where $C(I,S)=\p_{I,S}(u^{I,S})$ is a positive constant. Since $c(f)\in\R$ for any $\tilde{f}$, according to Lemma \ref{alem-2}, we have $c_{m}^{I,S}$=0.

The lemma is proved by induction on $m$.
\end{prfn}

\subsection{Quasi-local multi-vectors}

Let $\F$ be the quotient space $\A/\p\A$ whose  elements are called local functionals.
The coset of $f\in\A$ in $\F$ is denoted by $\int f\,dx$, and $f$ is called a density of this local functional.
The space $\F$ has a gradation induced from $\A$ since $\p$ is homogeneous of degree one.
Note that the operators $\delta_i$, $E$ satisfy $\delta_i\p=0$, $E\p=0$, so they induce maps from $\F$ to $\A$
which we still denote by $\delta_i$ and $E$.

We denote $\V^p=\mathrm{Alt}^p(\F,\F)$ where $\mathrm{Alt}^p(V_1, V_2)$ stands for
the linear space of all alternating $p$-linear maps from a linear space $V_1$ to another linear space $V_2$.
The elements of $\V^p$ are called the generalized $p$-vectors
over $\F$. We also use the notations $\V^0=\F$, $\V^{<0}=0$, and $\V=\bigoplus_{p\ge0}\V^p$.
The quasi-local multi-vectors that we are going to define below are certain kind of generalized multi-vectors
over $\F$ which are generated by the maps $E$ and $\delta_i$.

\begin{dfn}
Let $P\in\V^p$ be a generalized $p$-vector, we say that $P$ is quasi-local if the action of $P$ on $F_1, \dots, F_p\in\F$
takes the following form
\begin{align*}
&P(F_1, \dots, F_p)=\int \left(Q^{i_1 \dots i_p}_{s_1 \dots s_p}\,\p^{s_1}\delta_{i_1}(F_1)\dots\p^{s_p}\delta_{i_p}(F_p)\right.\\
&\quad +R^{j_1 \dots j_{p-1}}_{t_1 \dots t_{p-1}}\sum_{k=1}^p\left.(-1)^{k-1}E(F_k)\,\p^{t_1}\delta_{j_1}(F_1)\dots\hat{F}_k\dots
\p^{t_{p-1}}\delta_{j_{p-1}}(F_p)\right)\,dx,
\end{align*}
where $Q^{i_1 \dots i_p}_{s_1 \dots s_p},\ R^{j_1 \dots j_{p-1}}_{t_1 \dots t_{p-1}}\in\A$. $P$ is called
local if the second term of the above expression does not appear.
\end{dfn}

In order to give a more simple and concrete description of the class of quasi-local multi-vectors,
we introduce a family of super variables $\theta_i^s, \zeta$, where $i=1, \dots, n$, $s=0, 1, 2, \dots$, and define
\[\T=\A \otimes \wedge^*(V), \mbox{ where } V=\bigoplus_{i,s}\left(\mathbb{R}\theta_i^s\right)\oplus\mathbb{R}\zeta.\]
The gradation on $\A$ can be extended to a gradation on $\T$ as follow:
\begin{equation}
\deg u^{i,s}=\deg \theta_i^s=s,\ \deg f(u)=\deg \theta_i=\deg\zeta=0,\label{zh-bj-1}
\end{equation}
we still denote the completion of $\T$ w.r.t. the above gradation by $\T$.

The natural gradation of the exterior algebra $\wedge^*(V)$ yields another gradation on $\T$ as follows:
\[\T=\bigoplus_{p\ge0}\T^p,\ \T^p=\A\otimes\wedge^p(V).\]
In other words, elements of $\T^p$ have the form
\[P=P^{i_1 \dots i_p}_{s_1 \dots s_p}\theta_{i_1}^{s_1}\dots\theta_{i_p}^{s_p}+
X^{j_1 \dots j_{p-1}}_{t_1 \dots t_{p-1}}\,\zeta\,\theta_{j_1}^{t_1}\dots\theta_{j_{p-1}}^{t_{p-1}},\]
where $P^{i_1 \dots i_p}_{s_1 \dots s_p}, X^{j_1 \dots j_{p-1}}_{t_1 \dots t_{p-1}}\in\A$,
these coefficients are uniquely determined by $P$ if we require that they are antisymmetric w.r.t. the exchange of the indices $(i_k, s_k)\leftrightarrow (i_l,s_l)$
and $(j_k, t_k)\leftrightarrow (j_l,t_l)$.
An element $P\in\T$ is called local if it does not depend on $\zeta$.
The subspace of $\T$ consists of local elements is denoted by $\T_{loc}$

\begin{dfn}
We introduce a derivation $\hat{\p}:\T\to\T$
\begin{equation}
\hat{\p}=u^{i,s+1}\p_{i,s}+\theta_i^{s+1}\p^i_s-\left(u^{i,1}\theta_i\right)\p_{\zeta},
\end{equation}
and denote by $\E$ the quotient space $\T/\hat{\p}\T$.
Here $\p_{i,s}, \p^i_s, \p_\zeta$ are derivations on $\T$ given  by
\[\p_{i,s}=\frac{\p}{\p u^{i,s}},\quad \p^i_s=\frac{\p}{\p \theta_i^s},\quad \p_\zeta=\frac{\p}{\p \zeta}.\]
The coset of $P\in\T$ in $\E$ is denoted by $\int P\,dx$, and $P$ is called a density of $\int P\,dx$.
A coset $\int P\,dx \in\E$ is called local if it possesses a local density $P\in\T_{loc}$.
We denote the subspace of $\E$ consists of local elements by $\E_{loc}$.
Since $\hat{\p}$ is homogeneous w.r.t. both gradations of $\T$, the space $\E$ also has two gradations.
\end{dfn}
Note that the restriction of $\hat{\p}$ on $\A=\T^0$ coincides with $\p$, so it is not necessary to distinguish $\p$ and $\hat{\p}$.
We will denote $\hat{\p}$ by $\p$ henceforth.

The main result of the present subsection is the following theorem.

\begin{thm}\label{thm-jpp}
The space $\E$ is isomorphic to the space of quasi-local multi-vectors.
\end{thm}
\begin{prf}
We first define a map $\jmath:\T^p\to\mathrm{Alt}^p(\A,\A)$ as follow:
\begin{align}
&\jmath(P)(f_1, \dots, f_p)=\p^{i_p}_{s_p}\dots\p^{i_1}_{s_1}(P)\cdot\p^{s_1}\delta_{i_1}(f_1)\dots\p^{s_p}\delta_{i_p}(f_p)\nn\\
&+\p^{j_{p-1}}_{t_{p-1}}\dots\p^{j_1}_{t_1}\p_{\zeta}(P)\cdot\sum_{k=1}^p(-1)^{k-1}E(f_k)\,
\p^{t_1}\delta_{j_1}(f_1)\dots\hat{f}_k\dots\p^{t_{p-1}}\delta_{j_{p-1}}(f_p).\label{eq-jpp}
\end{align}
The identities $\p^i_s\p^j_t+\p^j_t\p^i_s=0$, $\p^i_s\p_\zeta+\p_\zeta\p^i_s=0$ imply that
$\jmath(P)$ is alternating.
Since $\delta_i \p=0$, $E \p=0$, the map $\jmath$ induces a map $\jmath':\T^p\to\mathrm{Alt}^p(\F,\A)$.

The desired isomorphism is defined as $\jmath'':\E^p \to \V^p$
\begin{equation}\label{zh-09-02}
\jmath''\left(\int P\,dx\right)(F_1, \dots, F_p)=\int \jmath'(P)(F_1, \dots, F_p)\,dx,
\end{equation}
This definition is independent of the choice of the density $P$, since we have
\[\jmath(\p Q)(f_1, \dots, f_p)=\p\left(\jmath(Q)(f_1, \dots, f_p)\right),\]
where $Q\in\T$ and $f_1, \dots, f_p\in\A$.

It is easy to see that the image of $\jmath''$ is just the subspace of quasi-local multi-vectors, the nontrivial
part of the theorem is to prove that $\jmath''$ is injective,
which is equivalent to the following lemma:
\begin{lem}\label{lem-hd}
If for $P\in\T^p$ and any $f_1, \dots, f_p\in\A$ we have
\[\jmath(P)(f_1, \dots, f_p)\sim 0,\]
then there exists $Q\in\T$ such that $P=\p Q$.
\end{lem}
The proof of the lemma will be given below, the theorem is then proved.
\end{prf}

In order to prove Lemma \ref{lem-hd}, we first need to prove the following two lemmas.

\begin{lem}\label{lem-4}
Lemma \ref{lem-hd} holds true when $p=0, 1, 2$.
\end{lem}
\begin{prf}
The case with $p=0$ is trivial. We begin with $p=1$.
Let $P\in\T^1$, so
\[P=Y^i_s\theta_i^s+h\,\zeta \sim X^i\theta_i+h\,\zeta,\]
where $X^i=\sum_{s\ge0}(-\p)^s Y^i_s$. Let $f\in\A$, then
\[\jmath(P)(f)\sim X^i \delta_i(f)+h\,E(f)\sim0.\]
Take $f=1$, we obtain $h\sim0$, so there exists $h'\in\A$ such that $h=\p h'$, then
\[P\sim X^i\theta_i+h\,\zeta\sim Z^i\theta_i,\]
where $Z^i=X^i+h'\,u^{i,1}$, and
\[\jmath(P)(f)\sim Z^i \delta_i(f)\sim0.\]
According to Lemma \ref{lem-2}, $Z^i=c\,u^{i,1}$, so
\[P\sim c\, u^{i,1}\theta_i =\p(-c\,\zeta)\sim0.\]
The $p=1$ case is proved.

Now let $P\in\T^2$, by performing integration by parts we have
\[P\sim P^{ij}_s \theta_i^s\theta_j+X^i\,\zeta\,\theta_i,\]
where $P^{ij}_s$ are uniquely determined by the antisymmetric condition
\[P^{ij}_s+\sum_{t\ge s}(-1)^t\binom{t}{s}\p^{t-s}P^{ji}_t=0.\]
Then
\[\jmath(P)(f,g)\sim 2P^{ij}_s\,\p^s\delta_i(f)\,\delta_j(g)+X^i\left(E(f)\,\delta_i(g)-\delta_i(f)\,E(g)\right)\sim0.\]
When $f=1$ we have $X^i\delta_i(g)\sim0$, so $X^i=c\,u^{i,1}$, thus
\[\jmath(P)(f,g)\sim \left(2P^{ij}_s\,\p^s\delta_i(f)+2\,c\,u^{j,1}E(f)\right)\,\delta_j(g)\sim0,\]
so there exists a differential operator $\tilde{c}:\A\to\R$ such that
\[2P^{ij}_s\,\p^s\delta_i(f)+2\,c\,u^{j,1}E(f)=\tilde{c}(f)u^{j,1},\]
but Lemma \ref{lem-3} tells us such $\tilde{c}$ must vanish, so
\[2P^{ij}_s\,\p^s\delta_i(f)+2\,c\,u^{j,1}E(f)=0,\]
then choose $f$ to be the monomials $1, u^i, u^iu^j, \dots$, one can obtain $P^{ij}_s=0$ and $c=0$, so $P\sim 0$.
The lemma is proved.
\end{prf}

\begin{lem}\label{lem-5}
Let $p\ge2$, $P\in\T^p_{loc}$, if for any $f_1, \dots, f_p\in\A$,
\[\jmath(P)(f_1, \dots, f_p)\sim0,\]
then there exists $Q\in\T^p_{loc}$ such that $P=\p Q$.
\end{lem}
\begin{prf}
Note that $P$ is homogeneous function of $\theta_i^s$, so
\begin{equation}\label{zh-10-10}
P=\frac1p \theta_i^s\p^i_s(P)\sim\frac1p\theta_i\delta^i(P)=\frac1pN(P),
\end{equation}
where 
\[\delta^i=\sum_{s\ge0}(-\p)^s\p^i_s,\] 
and 
\begin{equation}\label{zh-10-10-2}
N=\theta_i\delta^i
\end{equation} 
is the normalizing operator used in \cite{Ba}.

By the definitions of $\jmath$ and $\delta^i$, it is easy to see that
\[\jmath(P)(f_1, \dots, f_p)\sim\jmath(\delta^i(P))(f_2, \dots, f_p)\delta_i(f_1)\sim0,\]
so we have
\[\jmath(\delta^i(P))(f_2, \dots, f_p)=c(f_2, \dots, f_p)u^{i,1},\]
where $c(f_2, \dots, f_p)\in\R$ for any $f_2, \dots, f_p$.
According to Lemma \ref{lem-3}, the map $\jmath(\delta^i(P))$ is zero (here we used the condition $p\ge2$).
By choosing suitable $f_2, \dots, f_p$, one can obtain $\delta^i(P)=0$, so from
\eqref{zh-10-10} it follows that $P\sim0$. The lemma is proved.
\end{prf}

\begin{rmk}
The above lemma is not correct when $p=1$. In fact, $P=u^{i,1}\theta_i$ satisfies the condition, but $P=-\p \zeta$, and $\zeta \notin \T^1_{loc}$.
This fact also explains why the super variable $\zeta$ is so important.
\end{rmk}

\begin{prfn}{Lemma \ref{lem-hd}}
Since the cases with $p=0, 1, 2$ have been proved, we assume $p\ge3$.
It is easy to see
\[\jmath(P)(1, f_2, \dots, f_p)=-\jmath(\p_{\zeta}(P))(f_2, \dots, f_p)\sim0,\]
the derivative $\p_{\zeta}(P)$ satisfies the condition of Lemma \ref{lem-5}, so we have
\[\p_{\zeta}(P)=\p Q, \mbox{ where } Q\in\T^p_{loc}.\]
Let $P_0=P-\zeta\,\p_{\zeta}(P)$, then
\[P=P_0+\zeta\,\p_{\zeta}(P)=P_0+\zeta\,\p(Q)\sim P_0+u^{i,1}\theta_i\,Q\in\T^p_{loc}.\]
By applying Lemma \ref{lem-5} again, we proved the lemma.
\end{prfn}

Due to Theorem \ref{thm-jpp}, we denote the space of quasi-local multi-vectors also by $\E$.

\subsection{The Schouten-Nijenhuis bracket}

In this subsection we define the Schouten-Nijenhuis bracket on the space $\E$ of quasi-local multi-vectors,
it turns out to be the restriction of the Nijenhuis-Richardson bracket
defined on the space $\V$ of generalized multi-vectors to $\E$.

The original Nijenhuis-Richardson bracket was defined in \cite{NR-1,
NR-2} for finite dimensional vector spaces. In the following theorem we adopt it to the
infinite dimensional space $\V$ over the space of local functionals $\F$.

\begin{thm}\label{thm-bra}
There exists a unique bilinear map $[\ \, ,\ ]: \V^p\times\V^q\to\V^{p+q-1}$ satisfying the following conditions:
\begin{align}
&[P, F](F_2, \dots, F_p)=P(F, F_2, \dots, F_p), \label{nsb-1}\\
&[P,Q]=(-1)^{pq}[Q,P], \label{nsb-2}\\
&[[P,Q],F]+(-1)^{qp}[[Q,F],P]+[[F,P],Q]=0 \label{nsb-30}
\end{align}
for any $P\in\V^p$, $Q\in\V^q$, $F,\, F_2,\, \dots,\, F_p\in\F$. It is called the Nijenhuis-Richardson bracket of the generalized multi-vectors over $\F$.
\end{thm}

\begin{prf}
Note that for $P\in\V^p$, $Q\in\V^q$, $F_1, F_2, \dots, F_{p+q-1}\in\F$,
\begin{align*}
&[P,Q](F_1, F_2, \dots, F_{p+q-1})=[[P,Q],F_1](F_2, \dots, F_{p+q-1})\\
=&-\left((-1)^{qp}[[Q, F_1], P]+[[F_1, P],Q]\right)(F_2, \dots, F_{p+q-1}),
\end{align*}
so the bracket on $\V^p\times\V^q$ is determined by the brackets on $\V^{p-1}\times\V^q$ and $\V^p\times\V^{q-1}$, thus the uniqueness can be proved by
induction on $p+q$ immediately.

To prove the existence, we recall the product $\ow:\V^p\times\V^q\to\V^{p+q-1}$ defined in \cite{NR-2}:
\[P \ow Q(F_1, \dots, F_{p+q-1})=\sum_{I\in S_{p,q}}(-1)^{|I|}P(Q(F_{i_1},\dots, F_{i_q}), F_{i_{q+1}}, \dots, F_{i_{p+q-1}}),\]
where $S_{p,q}$ is the subset of the symmetry group $S_{p+q-1}$
\[S_{p,q}=\left\{I=(i_1, \dots, i_{p+q-1})\in S_{p+q-1}\left|\begin{array}{c} i_1<i_2<\dots<i_q \\ i_{q+1}<i_{q+2}<\dots<i_{p+q-1} \end{array}\right.\right\},\]
and $|I|$ denotes the parity of the permutation $I$. Note that this product is neither commutative nor associative.

We define the desired bracket as
\begin{equation}
[P,Q]=(-1)^{(p+1)q}\,P \ow Q+(-1)^p\,Q \ow P. \label{newsb}
\end{equation}
It is easy to verify that \eqref{newsb} satisfies the conditions \eqref{nsb-1} and \eqref{nsb-2}.
So we only need to prove the validity of the condition \eqref{nsb-30}.

Let $P\in\V^p$, $Q\in\V^q$, $F_1\in\F$. Denote
\[\tilde{P}=[P, F_1],\quad \tilde{Q}=[Q, F_1].\]
By definition of the product $\ow$, we have
\begin{align}
&P \ow Q(F_1, F_2, \dots, F_{p+q-1})\nn\\
=&P \ow \tilde{Q}(F_2, \dots, F_{p+q-1})+(-1)^{q+1}\tilde{P} \ow Q(F_2, \dots, F_{p+q-1}). \label{pcq}
\end{align}
It yields the identity
\[[[P,Q],F_1]=-(-1)^p[P,\tilde{Q}]-[\tilde{P},Q],\]
which is equivalent to \eqref{nsb-30}. The theorem is proved.
\end{prf}

\begin{cor}\label{cor-jcb}
The Nijenhuis-Richardson bracket satisfies the following graded Jacobi identity
\begin{equation}
(-1)^{pr}[[P,Q],R]+(-1)^{qp}[[Q,R],P]+(-1)^{rq}[[R,P],Q]=0, \label{nsb-3}
\end{equation}
for any $P\in\V^p,\ Q\in\V^q,\ R\in\V^r$.
\end{cor}

\begin{prf}
We prove the corollary by induction on $p+q+r$. When $r=0$, it is just the condition \eqref{nsb-30},
when $r>0$, we assume that the condition \eqref{nsb-3} holds true for any $p',q',r'$ with $p'+q'+r'<p+q+r$.
Let $P\in\V^p$, $Q\in\V^q$, $R\in\V^r$, $F\in\F$, then
\begin{align*}
&(-1)^{pr}[[[P, Q], R], F]\\
=&(-1)^{pr+p+q}[[P,Q],\tilde{R}]+(-1)^{pr+p}[[P, \tilde{Q}], R]+(-1)^{pr}[[\tilde{P}, Q], R],
\end{align*}
where $\tilde{K}=[K,F],\ K=P, Q, R$, so we have
\[[(-1)^{pr}[[P,Q],R]+(-1)^{qp}[[Q,R],P]+(-1)^{rq}[[R,P],Q],F]=0.\]
Thus the corollary follows from the property \eqref{nsb-1}. 
\end{prf}

\begin{rmk}
The proof of the graded Jacobi identity of the Richardson-Nijenhuis bracket given in \cite{NR-2} requires
that the vector space $\F$ is finite dimensional, while in the proof given above this condition is not necessary.
\end{rmk}

It is a nontrivial fact that the space $\E$ of quasi-local multi-vectors
is closed under the operation of Richardson-Nijenhuis bracket. In what follows we first define a bracket
on $\E$ by using the super variable description of the space of quasi-local multi-vectors,
and we prove that this bracket satisfies the conditions of Theorem \ref{thm-bra}. Then
the uniqueness property of Theorem \ref{thm-bra} shows that this bracket coincides with
the restriction of the Nijenhuis-Richardson bracket to $\E$. To this end, let us first introduce some
notations.

\begin{dfn}
For $i=1, \dots, n$, and integers $\alpha, s\ge0$, we define the higher generalized momentum operator, the higher energy operator and the energy operator
from $\T$ to $\T$ as
\begin{align*}
p_{i,\alpha,s}=&\sum_{\beta\ge0}(-1)^{\beta} \binom{\beta+s}{s} \p^{\beta}\p_{i,\alpha+\beta+s},\\
p^i_{\alpha,s}=&\sum_{\beta\ge0}(-1)^{\beta} \binom{\beta+s}{s} \p^{\beta}\p^i_{\alpha+\beta+s},\\
E_s=&\sum_{\alpha\ge1} \left(u^{i,\alpha}p_{i,\alpha,s}+\theta_i^{\alpha}p^i_{\alpha,s}\right),\\
E=&E_0-1,
\end{align*}
and assume that $p_{i,\alpha,-1}=\p_{i,\alpha-1}$, $p^i_{\alpha, -1}=\p^i_{\alpha-1}$. In particular, we have
\[E_{-1}=\sum_{\alpha\ge1}\left(u^{i,\alpha}\p_{i,\alpha-1}+\theta_i^{\alpha}\p^i_{\alpha-1}\right)=\p+\left(u^{i,1}\theta_i\right)\p_{\zeta},\]
and the higher Euler operators from $\T$ to $\T$ is defined as
\begin{equation}
\delta_{i,s}=p_{i,0,s},\ \delta^i_s=p^i_{0,s},\ \delta_i=\delta_{i,0},\ \delta^i=\delta^i_0.
\end{equation}
There is another useful operator $\hat{E}:\T\to\T$, which is defined as
\begin{equation}
\hat{E}=E+N,
\end{equation}
where $N$ is the normalizing operator defined in \eqref{zh-10-10-2}.
\end{dfn}

\begin{thm}\label{thm-main}
Define the bilinear map $\br:\E^p\times\E^q\to\E^{p+q-1}$ by
\begin{align}
&[\int P\,dx, \int Q\,dx]\nn\\
=&\int\left(\delta^i(P)\delta_i(Q)+(-1)^p\delta_i(P)\delta^i(Q)+\p_{\zeta}(P)\hat{E}(Q)+(-1)^p\hat{E}(P)\p_{\zeta}(Q)\right)\,dx. \label{mybra}
\end{align}
It satisfies the condition \eqref{nsb-1}-\eqref{nsb-30}, so it coincides with the Nijenhuis-Richardson bracket defined on $\V$.
\end{thm}
We call the above  bracket \eqref{mybra}  the {\em{Schouten-Nijenhuis bracket}} among quasi-local multi-vectors.
\begin{cor}
The bracket \eqref{mybra} satisfies the graded Jacobi identity
\begin{equation}
(-1)^{pr}[[P,Q],R]+(-1)^{qp}[[Q,R],P]+(-1)^{rq}[[R,P],Q]=0, \label{myjcb}
\end{equation}
for any $P\in\E^p,\ Q\in\E^q,\ R\in\E^r$.
\end{cor}

\begin{rmk}
The Schouten-Nijenhuis bracket can be restricted on the subspace $\E_{loc}$, the resulting bracket among local multi-vectors is equivalent to the
brackets defined by Getzler \cite{Ge} and by Kersten, Krasil$\,{}'$shchik, Verbovetsky \cite{KKV}.
\end{rmk}

\begin{rmk}
If $P, Q$ are of degree zero, then they become pairs of multi-vectors on the manifold $M$, and the bracket \eqref{mybra} degenerates to
\begin{align*}
[P,Q]=&\p^i(P)\p_i(Q)+(-1)^p\p_i(P)\p^i(Q)\\
&\quad+(1-q)\p_{\zeta}(P)Q+(1-p)(-1)^pP\p_{\zeta}(Q).
\end{align*}
This formula has appeared in the study of classical Jacobi structures, see \cite{IM} for example.
\end{rmk}

\begin{prfn}{Theorem \ref{thm-main}}
We first show that the bracket is well-defined, i.e. it is independent of the choices of $P, Q\in\T$.
Introduce a bilinear map $\br_{pr}:\T^p\times\T^q\to\T^{p+q-1}$
\[[P,Q]_{pr}=\delta^i(P)\delta_i(Q)+(-1)^p\delta_i(P)\delta^i(Q)+\p_{\zeta}(P)\hat{E}(Q)+(-1)^p\hat{E}(P)\p_{\zeta}(Q),\]
then it is easy to see that $[P, Q]_{pr}\sim D_P(Q)$, where $D_P:\T\to\T$ is a first order differential operator
\begin{equation}\label{zh-09-1}
D_P=-\p_{\zeta}(P)+(-1)^p \hat{E}(P)\p_{\zeta}+X^{i,s}\p_{i,s}+Y^s_i\p^i_s,
\end{equation}
$X^{i,s}, Y_i^s$ are defined by
\[X^{i,s}=\p^s\left(X^i\right)+\sum_{t=1}^s\p^{s-t}\left(\p_{\zeta}(P)u^{i,t}\right),\
Y_i^s=\p^s\left(Y_i\right)+\sum_{t=1}^s\p^{s-t}\left(\p_{\zeta}(P)\theta_i^t\right),\]
and $X^i=\delta^i(P)$, $Y_i=(-1)^p\left(\delta_i(P)-\theta_i\p_{\zeta}(P)\right)$.

Note that the operator $D_P$ satisfies $D_P\p=\p\left(D_P-D_P(1)\right)$, so
\[[P, \p (Q)]_{pr}\sim D_P\p(Q)=\p\left(D_P-D_P(1)\right)(Q)\sim0,\]
which implies that the definition is independent of the choice of $Q$. Then by using the fact that
\[[P,Q]_{pr}=(-1)^{pq}[Q,P]_{pr},\]
it is easy to see that this definition is also independent of the choice of $P$, so it is
well-defined on $\E^p\times \E^q$.

We next verify that the bracket defined in \eqref{mybra} satisfies the conditions \eqref{nsb-1}-\eqref{nsb-30}. In fact, the conditions \eqref{nsb-1} and
\eqref{nsb-2} are easy to verify, so we omit their proof here.
The condition \eqref{nsb-30} is a consequence of the following lemma:
\begin{lem} \label{DDf}
Let $P\in\T^p$, $f\in\A$, then we have the following identity:
\begin{equation} \label{zpf}
Z_{P,f}=D_{D_f(P)}+D_f\,D_P+(-1)^pD_P\,D_f=0.
\end{equation}
\end{lem}
The theorem then follows from this lemma, which will be proved below.
\end{prfn}

In order to prove Lemma \ref{DDf}, we first need to prove some other lemmas.

\begin{lem}\label{lem-9}
For $X,Y\in\E^1$, $F\in\F$, we have
\[[X,Y](F)=X(Y(F))-Y(X(F)),\]
here $X(\cdot)$ means $\jmath''(X)(\cdot)$, we omit $\jmath''$ from now on.
\end{lem}
\begin{prf}
The general form of an element in $\E^1$ reads
\[X=\int\left(X^i\,\theta_i+a\,\zeta\right)\,dx,\]
so we need to prove the lemma for the following three cases:
\begin{enumerate}
\item $X=\int\left(X^i\,\theta_i\right)\,dx$, $Y=\int\left(Y^j\,\theta_j\right)\,dx$;
\item $X=\int\left(X^i\,\theta_i\right)\,dx$, $Y=\int\left(b\,\zeta\right)\,dx$;
\item $X=\int\left(a\,\zeta\right)\,dx$, $Y=\int\left(b\,\zeta\right)\,dx$;
\end{enumerate}
We give below the proof for the third case, the proofs for the other two cases are similar.

According to the definition of $\br$, we have
\begin{align*}
&[X,Y]=\int\left(a\,\hat{E}(b\,\zeta)-b\,\hat{E}(a\,\zeta)\right)\,dx\\
=&\int\sum_{s\ge0}\sum_{t\ge1}(-1)^s\left(a\,u^{i,t}\p^s\p_{i,s+t}\left(b\,\zeta\right)-b\,u^{i,t}\p^s\p_{i,s+t}\left(a\,\zeta\right)\right)\,dx\\
=&\int\sum_{s\ge0}\sum_{t\ge1}\left(\p^s\left(a\,u^{i,t}\right)b_{(i,s+t)}-\p^s\left(b\,u^{i,t}\right)a_{(i,s+t)}\right)\,\zeta\,dx\\
=&\int\sum_{t\ge0}\left(\p^t\left(a\,E_t(b)\right)-\p^t\left(b\,E_t(a)\right)\right)\,\zeta\,dx.
\end{align*}
Here the last equality follows from the last identity of Lemma \ref{lem-idt1}.
In what follows we will also use this lemma.
So the left hand side reads (we omit the integral symbol to save space below)
\[[X,Y](F)\sim\sum_{t\ge0}\left(\p^t\left(a\,E_t(b)\right)-\p^t\left(b\,E_t(a)\right)\right)\,E(F).\]

On the other hand, let $F=\int\,f\,dx$, we have
\begin{align*}
&(X \,Y-Y\,X)(F)\sim a\,E(b\,E(f))-b\,E(a\,E(f))\\
=&\sum_{t\ge0}a(-1)^t\left(E_t(b)\p^t E(f)+\p^t(b)\,E_t E(f)\right)-a\,b\,E(f)\\
&\quad-\sum_{t\ge0}b(-1)^t\left(E_t(a)\p^t E(f)+\p^t(a)\,E_t E(f)\right)+a\,b\,E(f)\\
\sim&\sum_{t\ge0}\left(\p^t\left(a\,E_t(b)\right)-\p^t\left(b\,E_t(a)\right)\right)\,E(F)\\
&\quad+b\left(\sum_{t\ge0}\p^t\left(a\,E_t\,E(f)\right)-\sum_{t\ge0}(-1)^t\p^t(a)\,E_t\,E(f)\right),
\end{align*}
so we only need to prove the following identity
\[\sum_{t\ge0}\p^t\left(a\,E_t\,E(f)\right)=\sum_{t\ge0}(-1)^t\p^t(a)\,E_t\,E(f),\]
which is an easy corollary of the identity iii) in Lemma \ref{lem-idt3}. The lemma is proved.
\end{prf}

The following lemma is a generalization of Lemma \ref{lem-7}.
\begin{lem}\label{lem-7g}
Let $f\in\T$, if for any $g\in\A$ we have $f\,g\sim0$, then $f=0$.
\end{lem}
\begin{prf}
We first prove the lemma when $f$ is local. Without loss of generality, we assume that
$f\in\T^p_{loc}$ and is homogeneous
w.r.t. the gradation of $\T^p$ defined in \eqref{zh-bj-1}. The condition of the lemma can be written as
\begin{equation}\label{zh-bj-3}
g\,f=\p(h).
\end{equation}
The $p=0$ case is just Lemma \ref{lem-7}. Now we assume $p=1$, so that
\[f=\sum_{s=0}^m f^i_s\theta_i^s,\quad h=\sum_{s=0}^{m-1} h^i_s\theta_i^s+a\,\zeta.\]
Here $h$ is also homogeneous w.r.t. to both gradation of $\T$.
Compare the coefficients of $\theta_i^s$ and $\zeta$ on both sides of \eqref{zh-bj-3} we obtain
\begin{align}\label{zh-bj-2}
\p(a)=0,\ g\,f^i_0=\p\left(h^i_0\right)-a\,u^{i,1},\ g\,f^i_s=h^i_{s-1}+\p\left(h^i_s\right)\ (s\ge1),
\end{align}
from which we obtain $a\in\R$ and $g\,f^i_0\sim0$. By applying Lemma \ref{lem-7}, we also see that $f^i_0=0$ and $h^i_0=a u^i$.
From \eqref{zh-bj-2} it also follows that
\[a=a(g)=\frac1{u^i}\left( \sum_{l=0}^m (-\p)^{l-1} f^i_{l+1})\right) g,\]
so by using Lemma \ref{lem-3} we see that $a=0$ and thus $h^i_0=0$.
Now  by using the third equation of \eqref{zh-bj-2}, Lemma \ref{lem-7} and the homogeneity of
$h$ we arrive at $f=0$.
Thus  we proved the case when $p=1$.

When $p\ge2$, let $f\,g=\p(h)$ for certain $h\in\T^2$. We claim that $h\in\T^p_{loc}$.
In fact if $h=h_0+\zeta\,h_1$ with $h_0, h_1\in\T_{loc}$,
then $\p(h)=f\,g\in\T_{loc}$, so $\p(h_1)=0$ which implies that $h_1$ is a constant, and thus  $\zeta\,h_1\in\T^1$. Since $h\in \T^2$ we must have $h_1=0$.
Now for the local multi-vector $h\in\T^p_{loc}$ one can prove that $\delta^i\p(h)=0$, so
\begin{equation}\label{dfg}
0=\delta^i\left(f\,g\right)=\sum_{t\ge0}(-1)^t\,\delta^i_t(f)\p^t(g).
\end{equation}
Since we have  assumed that  $f$ is homogeneous in both gradations of $\T$,  there are only finite number
of nonzero terms in the right hand side of the above equation.
We regard the equation \eqref{dfg} as a system of linear equations of $(-1)^t\,\delta^i_t(f)$, then by choosing linear independent $g$'s and using
the property of Wronskian determinants, we obtain $\delta^i_t(f)=0$ for all $t\ge0$. Then the following identities
\[f=\frac1p\theta_i^s\p^i_s(f),\ \p^i_t(f)=\sum_{s\ge0}\binom{s+t}{s}\p^s\delta^i_{s+t}(f)\]
imply that $f=0$. Thus we proved the case when $f$ is local.

For the nonlocal case we note that $[\p,\p_{\zeta}]=0$, so
\[f \, g\sim0 \quad \Rightarrow \quad \p_{\zeta}(f)\,g\sim0,\]
by applying the above result to the local $\p_\zeta f$
we obtain  $\p_{\zeta}(f)=0$, so $f\in\T_{loc}$. By applying the above result again
to $f$ we finish the proof of the lemma.\end{prf}

\begin{rmk}
The identity $\delta^i\p(h)=0$ is not true if $\p_{\zeta}h\ne0$.
In general, the identities listed in Lemma \ref{lem-idt1} should be modified when the operators
act on $\T$.
\end{rmk}

\begin{prfn}{Lemma \ref{DDf}}
It is easy to see that the operator $Z_{P,f}$ is a first order differential operator and satisfies
\[Z_{P,f}\p=\p\left(Z_{P,f}-Z_{P,f}(1)\right),\]
so we only need to prove the validity of the following four identities
\[Z_{P,f}(1)=0,\ Z_{P,f}(\zeta)=0,\ Z_{P,f}(u^i)=0,\ Z_{P,f}(\theta_i)=0.\]

The first one is easy to prove. In fact, we have
\[Z_{P,f}(1)=-\left(\p_{\zeta}\,D_f+D_f\,\p_{\zeta}\right)(P)=0,\]
here we used the fact that $D_Q(1)=-\p_{\zeta}(Q)$ and $D_f=E(f)\p_{\zeta}+\p^s\left(\delta_i(f)\right)\p^i_s$.

For the third one we have
\[\left(Z_{P,f}-Z_{P,f}(1)\right)(u^i)=\left(\delta^i\,D_f+D_f\,\delta^i\right)(P).\]
The last expression equals zero since by using the fact that  $[D_f, \p]=0$  we obtain
\[\delta^i\,D_f+D_f\,\delta^i=\sum_{s\ge0}(-\p)^s\left(\p^i_s\,D_f+D_f\,\p^i_s\right)=0.\]

For the other two identities, we first consider the case when $p=1$.
Let $P\in\T^1$ and $Y=Y^i\theta_i+h\,\zeta\in\T^1$, Lemma \ref{lem-9} implies
\[0=\int Z_{P,f}(Y)\,dx=\int\left(Z_{P,f}(\zeta)h+Z_{P,f}(\theta_i)Y^i\right)\,dx,\]
then Lemma \ref{lem-7} implies that $Z_{P,f}(\zeta)=0$ and $Z_{P,f}(\theta_i)=0$, so we have $Z_{P,f}=0$.

Finally for $P\in\T^p$ and $Y=Y^i\theta_i+h\,\zeta\in\T^1$ we have
\[0=\int (-1)^p\,Z_{Y,f}(P)\,dx=\int Z_{P,f}(Y)\,dx.\]
By using the same argument as we used for the $p=1$ case and by using Lemma \ref{lem-7g} we obtain $Z_{P,f}(\zeta)=0$ and $Z_{P,f}(\theta_i)=0$, so we have $Z_{P,f}=0$.

The lemma is proved.
\end{prfn}

\subsection{Miura type transformations}

On the differential algebra $(\A, \p)$ there is a natural class of coordinate transformations which are called Miura type transformations in \cite{DZ},
they induce transformations of the graded Lie algebra $(\E, \br)$.  We are to
give some useful formulae for these transformations in this subsection.

\begin{dfn}
Let $\bar{u}^1, \dots, \bar{u}^n\in\A$ be $n$ differential polynomials, and $\bar{u}^i_0$ be the degree zero components of $\bar{u}^i$, if
\[\det\left(\frac{\p \bar{u}^i_0}{\p u^j}\right)\ne0\]
then the map $(u^1, \dots, u^n) \mapsto (\bar{u}^1, \dots, \bar{u}^n)$ defines a coordinate transformation on $\A$ which is called a Miura type transformation.
When $\bar{u}^i=\bar{u}^i_0$, this transformation is called a Miura type transformation of the first kind; when $\bar{u}^i_0=u^i$, this transformation
is called a Miura type transformation of the second kind.
\end{dfn}

It is easy to see that every Miura type transformation is the composition of a first kind Miura type transformation and a second kind one.
A first kind Miura type transformation is just a change of local coordinates on the base manifold $M$, and a general Miura type transformation is an analogue
of change of local coordinates on the jet space $J^{\infty}(M)$, note that $\A$ is similar to but different from $C^\infty(J^\infty(M))$.

Given a Miura type transformation $(\bar{u}^i)$, we can define $\bar{u}^{i,s}=\p^s\left(\bar{u}^i\right)$, and express $u^i$ as differential polynomials in
$(\bar{u}^{i,s})$, then every element of $\A$ can be expressed as differential polynomial in $(\bar{u}^{i,s})$, so we can define the following derivation
\[\bar{\p}_{i,s}=\frac{\p}{\p \bar{u}^{i,s}}:\A\to\A,\]
and the derivation $\p:\A\to\A$ now reads
\[\p=\sum_{s\ge0}\bar{u}^{i,s+1}\bar{\p}_{i,s}.\]
Furthermore, one can define operators $\bar{p}_{i,\alpha,s}$, $\bar{E}_s$ and $\bar{E}$ by replacing $\p_{i,s}$, $u^{i,s}$ by $\bar{\p}_{i,s}$,
$\bar{u}^{i,s}$ respectively in the original definition of $p_{i,\alpha, s}$, $E_s$, $E$ given in Section \ref{sec-2}.

\begin{lem}\label{lem-mde}
Let $f\in\A$ be a differential polynomial, then the following identities hold true:
\begin{align}
\delta_i(f)=&\sum_{s\ge0}(-\p)^s\left(\p_{i,s}(\bar{u}^j)\,\bar{\delta}_j(f)\right),\label{midt-1}\\
E(f)=&\bar{E}(f)-\sum_{s\ge0}\sum_{\alpha\ge1}u^{i,\alpha}(-\p)^s \left(\p_{i,s+\alpha}(\bar{u}^j)\,\bar{\delta}_j(f)\right).\label{midt-2}
\end{align}
\end{lem}
\begin{prf}
We denote the right hand sides of \eqref{midt-1} and \eqref{midt-2} by $g_i$ and $h$ respectively. According to Corollary \ref{lem-8}, we only need to prove that
for any first order differential operator $D:\A\to\A$ satisfying $D\p=\p\left(D-D(1)\right)$, one has
\begin{equation}
D(f)\sim X^i\,g_i-a\,h, \label{eq-df}
\end{equation}
where $a=D(1)$ and $X^i=\left(D-D(1)\right)(u^i)$.

On the other hand, if we apply Lemma \ref{lem-6} in the coordinate system $(\bar{u}^{i,s})$, we obtain
\begin{equation}
D(f)\sim \bar{X}^j\,\bar{\delta}_j(f)-\bar{a}\,\bar{E}(f), \label{eq-dfp}
\end{equation}
where $\bar{a}=D(1)=a$, and
\[\bar{X}^j=\left(D-D(1)\right)(\bar{u}^j)=\sum_{s\ge0}\left(\p^s(X^i)+\sum_{t=1}^s\p^{s-t}\left(a\,u^{i,t}\right)\right)\p_{i,s}(\bar{u}^j).\]
Then one can obtain \eqref{eq-df} from \eqref{eq-dfp} via integration by parts.
\end{prf}

\begin{emp}
Let $(\bar{u}^i)$ be a Miura type transformation of the first kind, then the identities \eqref{midt-1} and \eqref{midt-2} read
\[\delta_i(f)=\p_{i,0}(\bar{u}^j)\,\bar{\delta}_j(f),\ E(f)=\bar{E}(f).\]
The first one gives the classical transformation formula of variational derivatives, while the second one means that the energy is independent of the choice
of coordinates system.
\end{emp}

If we regard $\delta_i$ and $E$ as super variables $\theta_i$ and $\zeta$, the above lemma gives the transformation rule of super variables under Miura type transformations.
More precisely, we introduce a new family of super variables $\bar{\theta}_i^s, \bar{\zeta}$, and define the space
\[\bar{\T}=\A\otimes\wedge^*(\bar{V}), \mbox{ where } \bar{V}=\bigoplus_{i,s}\left(\mathbb{R}\bar{\theta}_i^s\right)\oplus\mathbb{R}\bar{\zeta}.\]
We can also define the quotient space $\bar{\E}=\bar{\T}/\p\bar{\T}$, and prove that $\bar{\E}$ is isomorphic to the space of quasi-local multi-vectors in $\V$,
so $\bar{\E}$ is also isomorphic to $\E$.

According to Lemma \ref{lem-mde}, we introduce the following isomorphism $\bar{\Phi}:\T\to\bar{\T}$ of super commutative algebras
\begin{align*}
\bar{\Phi}(f)=&f,\ f\in\A,\\
\bar{\Phi}(\theta_i^s)=&\p^s\sum_{t\ge0}(-\p)^t\left(\p_{i,t}(\bar{u}^j)\,\bar{\theta}_j\right),\\
\bar{\Phi}(\zeta)=&\bar{\zeta}-\sum_{s\ge0}\sum_{\alpha\ge1}u^{i,\alpha}(-\p)^s \left(\p_{i,s+\alpha}(\bar{u}^j)\,\bar{\theta}_j\right).
\end{align*}

\begin{thm}\label{thm-miura}
Let $P\in\V^p$ be a quasi-local $p$-vector, suppose the preimages of $P$ in $\E$ and $\bar{\E}$, under the isomorphism \eqref{zh-09-02}, are $\int \alpha\,dx$ and $\int \bar{\alpha}\,dx$ respectively,
then we have $\bar{\alpha}\sim \bar{\Phi}(\alpha)$.
\end{thm}
\begin{prf}
One can prove the theorem by acting $P$ on arbitrary $F_1, \dots, F_p\in\F$
and by using Lemma \ref{lem-mde}.
\end{prf}

For a Miura type transformation of the second kind, the formula given in Theorem \ref{thm-miura} is not ease to use.
The remaining part of this subsection is devoted to give another simpler formula for the transformation rule of the second kind Miura type transformations.
Let us first introduce the following notation:
\begin{dfn}\label{dfn-zh}
Let $F\in\T\mbox{ or }\E$, one can decompose $F$ w.r.t. the gradation \eqref{zh-bj-1} as follows:
\[F=F_k+F_{k+1}+F_{k+2}+\dots,\]
where $F_k\ne0$ and $\deg F_{l}=l$, then the integer $k$ is called the order of $F$, and is denoted by $\nu(F)$.
\end{dfn}

\begin{lem}
Let $(\bar{u}^i)$ be a second kind Miura type transformation, then there exists an $X\in\E^1_{loc}$ such that $\nu(X)>0$, and
\[\bar{u}^i=\exp(D_X)(u^i),\]
where $D_X$ is the first order differential operator associated to $X$ (see \eqref{zh-09-1} for the definition).
\end{lem}
\begin{prf}
Suppose $\bar{u}^i=u^i+F^i$, where $\nu(F^i)\ge k>0$, let
\[X_{(k)}=\int F^i_k\,\theta_i\,dx,\]
then we have
\[\exp\left(-D_{X_{(k)}}\right)(\bar{u}^i)=u^i+\tilde{F}^i,\]
where $\nu(\tilde{F}^i)\ge k+1>0$. So by induction on $k$, we can obtain a series
\[X_{(1)},\ X_{(2)},\ \dots \in\,\E^1_{loc},\]
such that
\[\bar{u}^i=\exp(D_{X_{(1)}})\exp(D_{X_{(2)}})\dots(u^i).\]
Finally, by using the Baker-Campbell-Hausdorff formula, and note that
\[\nu(X_{(1)})<\nu(X_{(2)})<\dots,\]
one can show that there exists $X\in\E^1_{loc}$ such that $\nu(X)>0$, and
\[\exp(D_X)=\exp(D_{X_{(1)}})\exp(D_{X_{(2)}})\dots.\]
The lemma is proved.
\end{prf}

From now on, we will denote a second kind Miura type transformation by $e^X$ for short, where $X\in\E^1_{loc}$ and satisfies $\nu(X)>0$.
When consider Miura type transformations of the second kind, it is convenient to identify $\E$ and $\bar{\E}$ and regard the induced isomorphism as automorphism on $\E$.

\begin{thm}\label{thm-miura2}
Let $e^X$ be a Miura type transformation of the second kind, then the induced automorphism of $\E$ is given by
\[\exp(-\ad_X):\E\to\E,\]
where $\ad_X=[X, \,\cdot\,]:\E\to\E$ is the adjoint action of $X$ on $\E$.
\end{thm}
\begin{prf}
Since both maps are given by exponential, we only need to compute the infinitesimal part of the transformation. For this we replace $X$ by $\e\,X$, then
the Miura type transformation reads
\[\bar{u}^i=u^i+\e\,X^i+O(\e^2),\]
where $X=\int X^i\,\theta_i\,dx$.

Let us consider the variation of the elements of $\cal{S}$ under Miura type transformations.
For $f\in\A$, we need to replace its independent variables $u^i$ by $\bar{u}^i-\e\,X^i+O(\e)^2$ and expand it into power series in $\e$, as the result we have
\[f \mapsto f+\e\,\Delta(f)+O(\e^2)=f-\e\left(\p^s(X^i)\p_{i,s}(f)\right)+O(\e^2).\]
For the super variables $\theta_i^s$ and $\zeta$, according to the definition of $\bar{\Phi}$, we have
\[\theta_i^s \mapsto \theta_i^s+\e\,\Delta(\theta_i^s)+O(\e^2),\ \zeta \mapsto \zeta+\e\,\Delta(\zeta)+O(\e^2),\]
where
\begin{align*}
&\Delta(\theta_i^s)=\p^s\sum_{t\ge0}(-\p)^t\left(\p_{i,t}(X^j)\theta_j\right),\\
&\Delta(\zeta)=-\sum_{s\ge0}\sum_{\alpha\ge1}u^{i,\alpha}(-\p)^s \left(\p_{i,s+\alpha}(X^j)\,\theta_j\right).
\end{align*}

Now let $P=\int \alpha\,dx\in\E$ be a quasi-local multi-vector, then the isomorphism $\bar{\Phi}$ becomes $P\mapsto P+\e\,\Delta(P)+O(\e)^2$, where
\[\Delta(P)\sim\Delta(u^{i,s})\p_{i,s}(\alpha)+\Delta(\theta_i^s)\p^i_s(\alpha)+\Delta(\zeta)\p_{\zeta}(\alpha),\]
then by using Theorem \ref{thm-main} and integration by parts, one can obtain
\[\Delta(P)=-[X, P].\]
The theorem is proved.
\end{prf}

\subsection{Reciprocal transformations}

In this subsection, we define the reciprocal transformation and give the transformation formula for quasi-local multi-vectors with the help of super variables.

Let $\rho\in\A$ be an invertible element, we define a new derivation
\begin{equation}
\tp=\rho^{-1}\p,
\end{equation}
and the new quotient space $\tF=\A/\tp\A$.
The coset of $\tilde{f}\in\A$ in $\tF$ is denoted by $\int \tilde{f}\,d\tx$.
It is easy to see that there is an isomorphism
\[\Phi_0:\F\to\tF,\ \int f\,dx\mapsto \int \rho^{-1}f\,d\tx.\]
We denote $\tV=\mathrm{Alt}^*(\tF, \tF)$, then one can extend
the isomorphism $\Phi_0:\F\to\tF$ to an isomorphism $\Phi:\V\to\tV$ such that
\begin{equation}
\Phi([P, Q])=[\Phi(P), \Phi(Q)]
\end{equation}
holds true for any $P, Q\in \V$.
In fact, for any $P\in\V^p$ the  action of $\Phi(P)$ on $\tilde{F}_1, \dots, \tilde{F}_p\in\tF$ is given as follows:
\[\Phi(P)(\tilde{F}_1, \dots, \tilde{F}_p)=\Phi_0\left(P\left(\Phi_0^{-1}\left(\tilde{F}_1\right), \dots, \Phi_0^{-1}\left(\tilde{F}_p\right)\right)\right).\]
By the definition of $\ow$, one can obtain $\Phi(P \ow\,Q)=\Phi(P)\ow\,\Phi(Q)$,
which implies that $\Phi([P, Q])=[\Phi(P), \Phi(Q)]$.

\begin{dfn}
The isomorphism $\Phi:\V\to\tV$ is called the reciprocal transformation w.r.t. $\rho$.
Let $\rho_0$ be the degree zero component of $\rho$, if $\rho=\rho_0$ then $\Phi$ is called a reciprocal transformation of the first kind, and if $\rho_0=1$
then $\Phi$ is called a reciprocal transformation of the second kind.
\end{dfn}

When restricted to the subspace $\E\subset\V$ the reciprocal transformation $\Phi$ has a simple computation formula, we will give it below after
some preparations.

We introduce a new coordinate system on $\A$
\[\tu^i=u^i,\ \tu^{i,s}=\tp^s(\tu^i).\]
It is easy to see that the transformation $\{u^{i,s}\}\to\{\tu^{i,s}\}$ is invertible, and elements of $\A$ can be written as differential polynomials in $\tu^{i,s}$.
By defining the new derivations $\tp_{i,s}=\frac{\p}{\p \tu^{i,s}}$, one can represent $\tp$ as
\[\tp=\sum_{s\ge0}\tu^{i,s+1}\tp_{i,s}.\]
Furthermore, one can define operators $\tilde{p}_{i,\alpha,s}$, $\tilde{E}_s$ and $\tilde{E}$ by replacing $\p$, $\p_{i,s}$, $u^{i,s}$ by $\tp$, $\tp_{i,s}$,
$\tu^{i,s}$ respectively in the original definition of $p_{i,\alpha, s}$, $E_s$, $E$ given
in Section \ref{sec-2}.

\begin{lem}\label{lem-dede}
Let $f\in\A$, $\tilde{f}=\rho^{-1}f$, then the following identities hold true
\begin{align}
\delta_i(f)=&\rho\,\tilde{\delta}_i(\tilde{f})-\sum_{s\ge0}(-\p)^s\left(\p_{i,s}(\rho)\tilde{E}(\tilde{f})\right),\label{ridt-1}\\
E(f)=&\rho\,\tilde{E}(\tilde{f})-\sum_{s\ge0}\sum_{\alpha\ge1}u^{i,\alpha}(-\p)^s\left(\p_{i,s+\alpha}(\rho)\tilde{E}(\tilde{f})\right). \label{ridt-2}
\end{align}
\end{lem}
\begin{prf}
We denote the right hand sides of \eqref{ridt-1} and \eqref{ridt-2} by $g_i$ and $h$. According to Corollary \ref{lem-8}, we only need to prove that
for any first order differential operator $D:\A\to\A$ satisfying $D\p=\p\left(D-D(1)\right)$, one has
\begin{equation}
D(f)\sim X^i\,g_i-a\,h, \label{eq-topr}
\end{equation}
where $a=D(1)$ and $X^i=\left(D-D(1)\right)(u^i)$. In fact, such $D$ corresponds to a quasi-local vector $X=X^i\theta_i-a\,\zeta\in\E^1$,
and $D(f)\sim X(F)$, where $F=\int f\,dx$.

Define $\tilde{F}=\Phi(F)=\int \tilde{f}\,d\tx$, then we have
\[[\Phi(X),\tilde{F}]=\Phi(X)(\tilde{F})=\Phi_0\left(X\left(\Phi_0^{-1}(\tilde{F})\right)\right)=\int \left(\rho^{-1}D\rho\right)(\tilde{f})\,d\tx.\]
Denote $\tilde{D}=\rho^{-1}D\rho$, one can check that $\tilde{D}$ is a first order differential operator satisfying
$\tilde{D}\tp=\tp\left(\tilde{D}-\tilde{D}(1)\right)$, then Lemma \ref{lem-6} implies
\[\tilde{D}(\tilde{f})\tilde{\sim}\tilde{X}^i\tilde{\delta}_i(\tilde{f})-\tilde{a}\,\tilde{E}(\tilde{f}),\]
where $\tilde{a}=\tilde{D}(1)=\rho^{-1}D(\rho)$, $\tilde{X}^i=\left(\tilde{D}-\tilde{D}(1)\right)(\tu^i)=X^i$, and
\[A\tilde{\sim}B \Leftrightarrow A-B\in\tp\A.\]

On the other hand,
\[D(f)\sim X(F)=\Phi_0^{-1}\left([\Phi(X),\tilde{F}]\right)=\Phi_0^{-1}\left(\int \tilde{D}(\tilde{f})\,d\tx\right)=\int \rho\,\tilde{D}(\tilde{f})\,dx,\]
so we have
\[D(f)\sim \rho\,X^i\tilde{\delta}_i(\tilde{f})-D(\rho)\,\tilde{E}(\tilde{f}),\]
then one can obtain \eqref{eq-topr} after integration by parts.
The lemma is proved.
\end{prf}

As we did in the last subsection for Miura type transformations, if we regard $\delta_i$ and $E$ as super variables $\theta_i$ and $\zeta$ the above lemma gives the transformation rule of super variables
under reciprocal transformations.
More precisely, we introduce a new family of super variables $\tilde{\theta}_i^s, \tilde{\zeta}$, and define the space
\[\tT=\A\otimes\wedge^*(\tilde{V}), \mbox{ where } \tilde{V}=\bigoplus_{i,s}\left(\mathbb{R}\tilde{\theta}_i^s\right)\oplus\mathbb{R}\tilde{\zeta}.\]
We can also define the quotient space $\tE=\tT/\tp\tT$ which is isomorphic to the space of quasi-local multi-vectors in $\tV$.

According to Lemma \ref{lem-dede}, we introduce the following isomorphism $\hat{\Phi}:\T\to\tT$ of super commutative algebras
\begin{align*}
\hat{\Phi}(f)=&f,\ f\in\A,\\
\hat{\Phi}(\theta_i^s)=&\left(\rho\tp\right)^s\left(\rho\,\tilde{\theta}_i-\sum_{t\ge0}(-\rho\tp)^t\left(\p_{i,t}(\rho)\tilde{\zeta}\right)\right),\\
\hat{\Phi}(\zeta)=&\rho\,\tilde{\zeta}-\sum_{s\ge0}\sum_{\alpha\ge1}u^{i,\alpha}(-\rho\tp)^s\left(\p_{i,s+\alpha}(\rho)\tilde{\zeta}\right).
\end{align*}

\begin{thm}\label{thm-reci}
The restriction of the reciprocal transformation $\Phi$ on $\E\subset\V$ is an isomorphism from $\E$ to $\tE$.
More precisely, let $\alpha \in \T^p$, we have
\begin{equation}
\Phi\left(\int \alpha\,dx\right)=\int \rho^{-1}\hat{\Phi}(\alpha)\,d\tx. \label{trsp}
\end{equation}
\end{thm}
\begin{prf}
The theorem can be easily proved by acting both sides of \eqref{trsp} on the local functionals $\tilde{F_1}, \dots, \tilde{F}_p\in\tF$
and using Lemma \ref{lem-dede}.
\end{prf}

\begin{emp}
We consider an evolutionary PDE
\begin{equation}\label{zh-09-03}
\p_tu^i=X^i, \ \mbox{where}\  X^i\in\A.
\end{equation}
It corresponds to a local vector $X=\int X^i\theta_i\,dx\in\E^1$.
Let $\rho\in\A$ be an invertible element, we can define the reciprocal transformation $\Phi$ w.r.t. $\rho$.
From Theorem \ref{thm-reci} it follows that
\begin{equation}
\Phi(X)=\int\left(X^i\tilde{\theta}_i-\rho^{-1}\p_t(\rho)\tilde{\zeta}\right)\,d\tx.
\end{equation}
If $\rho$ is a conserved density of $\p_t$, i.e. there exists $\sigma\in\A$ such that $\p_t(\rho)=\p(\sigma)$, then
\[\Phi(X)=\int\left(X^i-\sigma\,\tilde{u}^{i,1}\right)\tilde{\theta}_i\,d\tx.\]
So we obtain another evolutionary PDE
\begin{equation}
\p_{\tilde{t}}{\tilde{u}}^i=X^i-\sigma\,\tilde{u}^{i,1}.
\end{equation}
This equation coincides with the one that is obtained from \eqref{zh-09-03} by the following reciprocal transformation:
\[d\tx=\rho\,dx+\sigma\,dt,\ d\tilde{t}=dt.\]
\end{emp}

If $\Phi$ is a first kind reciprocal transformation, then the isomorphism $\hat{\Phi}$ is very simple:
\begin{equation}\label{zh-09-06}
\hat{\Phi}(f)=f,\ \hat{\Phi}(\theta_i^s)=\left(\rho\tp\right)^s\left(\rho\tilde{\theta}_i-\p_{i,0}(\rho)\tilde{\zeta}\right),\ \hat{\Phi}(\zeta)=\rho\tilde{\zeta}.
\end{equation}
When $\Phi$ is of the second class, we first identify $\E$ and $\tE$, and regard $\Phi$ as an automorphism on $\E$,
then we have the following formula.

\begin{thm}\label{thm-reci2}
Let $\rho=e^f$, where $f\in\A_{>0}$, and let $\Phi$ be the reciprocal transformation w.r.t. $\rho$, then for any $P\in\E$ we have
\[\Phi(P)=\exp(\ad_Y)(P),\]
where $Y=\int f\,\zeta\,dx\in\E^1$.
\end{thm}
\begin{prf}
The proof of this theorem is similar to that of Theorem \ref{thm-miura2}, so we omit the details here.
\end{prf}

\begin{rmk}
Theorem \ref{thm-miura2} and Theorem \ref{thm-reci2} show that both Miura type transformations of the second kind and
reciprocal transformations of the second kind
can be expressed as $\exp(\ad_X)$, where $X\in\E^1_{>0}$. Conversely, by using the Baker-Campbell-Hausdorff formula one can show that
every automorphism $\exp(\ad_X):\ \E\to \E$ with $X\in\E^1_{>0}$ can be represented as the composition of a
Miura type transformation of the second kind and a reciprocal transformation of the second kind.
This observation is very important when we study the classification problem of deformations of
Jacobi structures (see Section \ref{sec-dfm} below).
\end{rmk}

\section{Jacobi structures}

\subsection{Definition and examples} \label{sec-hydro}

A Jacobi structure on a finite dimensional manifold $M^n$, as it was introduced by Lichnerowicz in \cite{lich-2},
consists of a pair $(\Lambda, X)$ of a bivector $\Lambda$ and a vector field $X$ on $M$, they
satisfy the conditions
\begin{equation}
[\Lambda, \Lambda]=2 X\wedge \Lambda,\quad [X, \Lambda]=0.
\end{equation}
It is equivalent to a local Lie algebra structure defined on the space of smooth functions on $M$
via the following bracket:
\[
\{f, g\}=\Lambda(df, dg)+f X(g)-g X(f),\quad \forall f, g\in C^\infty(M).
\]
This bracket satisfies the Jacobi identity, and in general (when $X\ne 0$)
do not satisfies the Leibniz rule.
In this section, we generalize such structures to the infinite jet space of $M$, and in this way define Jacobi structures of evolutionary PDEs.

We denote the Lie algebra of derivations on $\A$ by $\Der(\A)$, and denote the centralizer of $\p$ in $\Der(\A)$
by $\Der'(\A)$.
It is easy to see that the space $\Der'(\A)$ consists of elements of the form
\[D=\sum_{s\ge0}\p^s\left(X^i\right)\p_{i,s},\quad X^i\in\A.\]
Each such a derivation corresponds to an evolutionary PDE with components in $\A$:
\begin{equation}
\p_tu^i=X^i,\ X^i\in\A. \label{epde}
\end{equation}
Note that
the action of $D$ on $f\in\A$ is just $\p_tf$.

Since elements of $\Der'(\A)$ commute with $\p$, $\F$ has a natural $\Der'(\A)$-module structure given by
\[\pi:\Der'(\A)\to\V^1=\mathrm{End}(\F, \F),\ \pi(D)\left(\int f\,dx\right)=\int D(f)\,dx.\]
The image of $\pi$ is just $\E^1_{loc}$ (since $D(f)\sim X^i\,\delta_i(f)$),
and it follows from Lemma \ref{lem-2} that  the kernel of $\pi$ is $\R\,\p$. So we have the following isomorphism:
\[\E^1_{loc}\cong \Der'(\A)/\R\,\p.\]
This quotient has a nice physical explanation. Given an evolutionary PDE \eqref{epde}, we can convert it to
\[\p_tu^i=X^i+c\,u^{i,1}\]
by performing the following Galilean transformation
\[x\mapsto x-c\,t,\ t\mapsto t,\]
so the space $\E^1_{loc}$ is just the space of equivalence classes of evolutionary PDEs modulo Galilean transformations.
In this paper, we consier these equivalence classes only, and call elements of $\E^1_{loc}$ evolutionary PDEs for short.

\begin{dfn}
A quasi-local bivector  $P\in\E^2$ is called a Jacobi structure if  $[P, P]=0$.
An evolutionary PDE $X\in\E^1_{loc}$ is said to possess a Jacobi structure
if there exists a Jacobi structure $P$ and a local functional $F\in\F$ such that $X=[P,F]$,
here  $F$ is called the Hamiltonian.
\end{dfn}

By using integration by parts one can always represent a Jacobi structure in the following form:
\[P=\int\left(\frac12\theta_i\left(\alpha^{ij}_s\p^s\right)\theta_j+\zeta\,X^i\theta_i\right) dx,\]
where $\alpha^{ij}_s,\ X^i\in\A$, and the matrix differential operator $\alpha=\left(\alpha^{ij}_s\p^s\right)$
satisfies the skew-symmetry condition
\begin{equation}
\alpha+\alpha^{\dagger}=0 \label{adagger}
\end{equation}
with $\alpha^{\dagger}=\left(((-\p)^s\alpha^{ji}_s)\right)$. The bivector $P_0=\frac12\theta_i\alpha^{ij}\theta_j\in\T^2_{loc}$ is called the local part of $P$
and $X=X^i\theta_i\in\T^1_{loc}$ is called, following the notation of \cite{FP},  the structure flow of $P$. In what follows a Jacobi structure $P$ is often
represented as $P=(P_0, X)$, or $P\sim P_0+\zeta\,X$.

\begin{lem}\label{lem-PX}
Let $X\in\E^1_{loc}$, $P\in\E^2$ be a Jacobi structure of $X$, then the structure flow of $P$ is a symmetry of $X$.
\end{lem}
\begin{prf}
By definition, there exists $F\in\F$ such that $X=[P, F]$, so we have $[X, P]=0$.
Note that the derivation $\p_{\zeta}:\T\to\T$ satisfies $[\p_{\zeta}, \p]=0$, so it induces a map $\p_{\zeta}:\E\to\E$.
By the action of $\p_{\zeta}$, one can obtain
\[0=\p_{\zeta}\left([X, P]\right)=[X, \p_{\zeta}(P)],\]
where $\p_{\zeta}(P)$ is just the structure flow of $P$. The lemma is proved.
\end{prf}

\begin{lem}\label{lem-XH}
Let $P=(P_0, X)$ be a Jacobi structure, $H\in\F$ be a local functional, then $[P, H]$ is local if and only if $[X, H]=0$.
\end{lem}
\begin{prf}
The locality of $[P,H]$ means that $\p_{\zeta}[P,H]=0$, so we have
\begin{align*}
0=&\p_{\zeta}[P,H]\sim\p_{\zeta}\left(\delta^i(P)\delta_i(H)+\p_{\zeta}(P)\hat{E}(H)\right)\\
=&\p_{\zeta}\left(\delta^i(\zeta\,X)\right)\delta_i(H)\sim-[X,H],
\end{align*}
the lemma is proved.
\end{prf}

The above lemma means that for a given Jacobi structure $P=(P_0, X)$, the admissible Hamiltonians must be conserved quantities of the structure flow $X$.
This is quite different from the local case.

\begin{emp}\label{zh-bj-5}
The simplest Jacobi structures are the ones with degree zero. Let $P=(P_0, X)\in\E^2$ be a Jacobi structure such that $\deg(P)=0$, then
$P_0$ and $X$ are just a bivector and a vector field on $M$ respectively. The condition $[P, P]=0$ implies
\[[P_0, P_0]+2\,X\,P_0=0,\ [P_0, X]=0,\]
so the pair $(P_0, -X)$ gives a classical Jacobi structure on $M$ \cite{lich-2}.

Let $P\sim\frac12\alpha^{ij}\theta_i\theta_j+\zeta\,X^i\,\theta_i$ be a Jacobi structure of degree zero, and $\rho$ be a nowhere zero smooth function on $M$.
One can define the reciprocal transformation $\Phi$ w.r.t. $\rho$. The isomorphism $\hat{\Phi}$ is given by \eqref{zh-09-06}, so we obtain the reciprocal
transformation of $P$
\[\Phi(P)\tilde{\sim}\frac12\tilde{\alpha}^{ij}\tilde{\theta}_i\tilde{\theta}_j+\tilde{\zeta}\,\tilde{X}^i\,\tilde{\theta}_i,\]
where $\tilde{\alpha}^{ij}=\rho\,\alpha^{ij}$, $\tilde{X}^i=\rho\,X^i+\alpha^{ij}\p_{j,0}(\rho)$. This is in fact a conformal change of classical
Jacobi structures \cite{GL, DLM}.
\end{emp}

The Jacobi structures of degree $1$ are more interesting for us. Let $P=(P_0, X)$ be a Jacobi structure where
\begin{equation}
P_0=\frac12\left(g^{ij}(u)\theta_i\theta_j^1+\Gamma^{ij}_k(u)u^{k,1}\theta_i\theta_j\right),\quad  X=V^i_k(u)u^{k,1}\theta_i, \label{hydro}
\end{equation}
and we assume that $\det\left(g^{ij}\right)\ne0$. The skew-symmetry condition \eqref{adagger} is equivalent to
\[g^{ij}=g^{ji},\ \Gamma^{ij}_k+\Gamma^{ji}_k=g^{ij}_{(k,0)}=\frac{\p g^{ij}}{\p u^k},\]
so $g^{ij}$ is a contravariant metric on $M$, we denote its inverse by $g_{ij}$, and define $\Gamma^j_{kl}=-g_{ki}\Gamma^{ij}_l$,
$V_{kj}=g_{ki}V^i_j$.

The following theorem was first proved by Ferapontov in \cite{fera-1}, here we give another proof to illustrate the usage of our formula
of the Schouten-Nijenhuis bracket \eqref{mybra}.

\begin{thm}\label{thm-fera}
A bivector $P=(P_0, X)$ of the form \eqref{hydro} is a Jacobi structure if and only if the following conditions are satisfied
\begin{align}
&\Gamma^j_{kl}=\Gamma^j_{lk},\ V_{kj}=V_{jk},\ \nabla_kV_{lj}=\nabla_lV_{kj},\nn\\
&R_{ijkl}=g_{ik}V_{jl}+g_{jl}V_{ik}-g_{jk}V_{il}-g_{il}V_{jk},\label{four-c}
\end{align}
where $\nabla$ and $R_{ijkl}$ are the Levi-Civita connection and Riemannian curvature tensor of $g_{ij}$ respectively.
\end{thm}
\begin{prf}
We need to check that the condition $[P, P]=0$ is equivalent to the conditions listed in the theorem. From \eqref{mybra} it follows that
\[\frac12[P, P]\sim\delta^i(P_0+\zeta\,X)\delta_i(P_0+\zeta\,X)+X\,\hat{E}(P_0+\zeta\,X).\]
By using the following identities,
\begin{align*}
&\delta^i(\zeta\,X)=-\zeta\,\delta^i(X),\ \delta_i(\zeta\,X)=\zeta\,\delta_i(X)-\left(u^{k,1}\theta_k\right)\left(V^l_i\theta_l\right),\\
&\hat{E}(\zeta\,X)=\zeta\,X,\ \hat{E}(P_0)=2P_0,
\end{align*}
one can obtain $\frac12[P,P]\sim A-\zeta\,B$, where
\begin{align*}
A=&\frac12[P_0, P_0]_{pr}-\delta^i(P_0)\left(u^{k,1}\theta_k\right)\left(V^l_i\theta_l\right)+2\,X\,P_0,\\
B=&[P_0,X]_{pr}-\delta^i(X)\left(u^{k,1}\theta_k\right)\left(V^l_i\theta_l\right).
\end{align*}

The condition $A-\zeta\,B\sim0$ implies $B\sim0$. Note that $B\in\T^2_{loc}$, so we have $\delta^i(B)=0$ (see the proof of Lemma \ref{lem-7g}).
We rewrite $B$ in the following form
\[B=E^{ij}\theta_i^1\theta_j^1+F^{ij}_ku^{k,1}\theta_i\theta_j^1+H^{ij}_{kl}u^{k,1}u^{l,1}\theta_i\theta_j,\]
then we have
\begin{align}
0&=\p^i_2\delta^j(B)=\p^i_2\left(\p^j_0-\p\p^j_1\right)(B)=\p^j_1\p^i_1(B)=E^{ij}-E^{ji}, \label{con-E}\\
0&=\p_{k,2}\delta^j(B)=\p_{k,2}\left(\p^j_0-\p\p^j_1\right)(B)=-\p_{k,1}\p^j_1(B)=F^{ij}_k\theta_i, \label{con-F}
\end{align}
so we have $\delta^i(B)=\left(H^{ij}_{kl}-H^{ji}_{kl}\right)u^{k,1}u^{l,1}\theta_j=0$, which implies
\begin{equation}
H^{ij}_{kl}-H^{ji}_{kl}+H^{ij}_{lk}-H^{ji}_{lk}=0. \label{con-H}
\end{equation}

The condtions \eqref{con-E}, \eqref{con-F}, \eqref{con-H} show that $B=0$, so we have $A\sim0$ and $\delta^i(A)=0$.
Then by straightforward computation, one can obtain $\Gamma^j_{kl}=\Gamma^j_{lk}$ from $\p^i_2\delta^j(A)=0$, and $V_{kj}=V_{jk}$ from \eqref{con-E},
and $\nabla_kV_{lj}=\nabla_lV_{kj}$ from \eqref{con-F}, and $R_{ijkl}$ from $\p_{k,2}\delta^j(A)=0$. So we proved that the four conditions are necessary.

The sufficiency part can be proved by a straightforward but lengthy computation, we omit it here.
\end{prf}

\begin{dfn}
A bivector $P\sim P_0+\zeta\,X\in\E^2$ is called a Jacobi structure of hydrodynamic type
if $(P_0, X)$ is given by \eqref{hydro} and $\det\left(g^{ij}\right)\ne0$.
\end{dfn}

In what follows, we will study reciprocal transformations of hydrodynamic Jacobi structures.
Note that reciprocal transformations preserve the Schouten-Nijenhuis bracket, so the reciprocal
transformation of a Jacobi structure is still a Jacobi structure.
The following fundamental formula was given by Ferapontov and Pavlov.
\begin{lem}[\cite{FP}]\label{lem-FP}
Let $P=(P_0, X)$ be a Jacobi structure of hydrodynamic type, and $\rho$ be a nowhere zero smooth function on $M$.
We denote the reciprocal transformation w.r.t. $\rho$ by $\Phi$ and assume $\Phi(P)=(\tilde{P}_0, \tilde{X})$, where
\[\tilde{P}_0=\frac12\left(\tilde{g}^{ij}(u)\tilde{\theta}_i\tilde{\theta}_j^1+\tilde{\Gamma}^{ij}_k(u)\tilde{u}^{k,1}\tilde{\theta}_i\tilde{\theta}_j\right),
\ \tilde{X}=\tilde{V}^i_k(u)\tilde{u}^{k,1}\tilde{\theta}_i,\]
then we have
\begin{equation}
\tilde{g}^{ij}=\rho^2\,g^{ij},\ \tilde{V}^i_k=\rho^2 V^i_k+\rho\,\nabla^i\nabla_k\rho-\frac12(\nabla^l\rho)(\nabla_l\rho)\delta^i_k. \label{eq-fera}
\end{equation}
Here $\nabla$ is the Levi-Civita connection of the metric $g_{ij}$.
\end{lem}
In fact the \eqref{eq-fera} can also be easily obtained from the relation \eqref{zh-09-06} and our general transformation formula \eqref{trsp}.

\begin{thm}\label{zh-jht}
Let $P=(P_0, X)$ be a Jacobi structure of hydrodynamic type, then there exists a reciprocal transformation converting $P$ to a local Jacobi structure, i.e.
a Hamiltonian structure.
\end{thm}
\begin{prf}
From the condition \eqref{four-c}, one can obtain the Ricci and scalar curvatures of $g_{ij}$
\[R_{ij}(g)=(n-2)V_{ij}+\left(V^k_k\right)g_{ij},\ R(g)=2(n-1)V^k_k\]
which imply that both the Weyl tensor and the Cotton tensor of $g_{ij}$ vanish, so $g_{ij}$ is a conformally flat metric.

Let $g_{ij}=\rho^{-2}\,\eta_{ij}$, where $\eta_{ij}$ is a flat metric, and $\rho^{-2}$ is the conformal factor.
Note that $\rho^{-1}$ is nowhere zero on $M$, so we can define a reciprocal transformation $\Phi$ w.r.t. $\rho^{-1}$.
Let $\Phi(P)=Q=(Q_0, Y)$, then it is easy to see that the metric of $Q_0$ is just the flat
metric $\eta_{ij}$.
So without loss of generality we can assume that $g_{ij}$ has been converted to a flat metric $\eta_{ij}$.

When $n\ge3$, one obtain
\[V_{ij}=\frac{R_{ij}(\eta)}{n-2}-\frac{g_{ij}R(\eta)}{2(n-1)(n-2)}=0,\]
so the theorem has been proved in this case.

Now we assume $n=2$. In this case, the tensor $V_{ij}$ is not determined by the curvatures, so it may not be zero.
We choose $u^1, u^2$ to be the flat coordinates of metric $\eta$, then $\eta_{ij}$ becomes a constant symmetric matrix.
There are two possible cases up to linear transformations
\[(\eta_{ij})=\left(\begin{array}{cc}1 & 0 \\ 0 & 1 \end{array}\right)\  \mbox{ or }\ (\eta_{ij})=\left(\begin{array}{cc} 0 & 1\\ 1 & 0 \end{array}\right).\]
We only consider the first case since the second one is similar and easier.

We denote the components of the tensor $V_{ij}$ by
\[a=V_{11},\ b=V_{12}=V_{21},\ c=V_{22},\]
and $u^1=x, u^2=y$, then the conditions \eqref{four-c} imply
\[a+c=0,\ a_y=b_x,\ b_y=c_x.\]
So if we introduce $z=x+y\,\sqrt{-1}$, $f=a-b\,\sqrt{-1}$, then $f$ is holomorphic in $z$.

On the other hand, if there exists $\rho$ such that $\Phi(P)$ is local, then this $\rho$ must satisfy the following equation which is implied by the formula
\eqref{eq-fera}:
\[\rho^2 V_{ij}+\rho\,\nabla_i\nabla_j\rho-\frac12(\nabla^l\rho)(\nabla_l\rho)\eta_{ij}=0.\]
The components of the above equation read
\begin{equation}
\frac{\rho_x^2+\rho_y^2}{2\,\rho^2}-\frac{\rho_{xx}}{\rho}=a,\ -\frac{\rho_{xy}}{\rho}=b,\ \frac{\rho_x^2+\rho_y^2}{2\,\rho^2}-\frac{\rho_{yy}}{\rho}=-a.
\label{eq-ric}
\end{equation}
We introduce a new function
\[\chi=-\frac{\rho_x}{\rho}+\frac{\rho_y}{\rho}\,\sqrt{-1},\]
then the equations \eqref{eq-ric} imply that $\chi$ is holomorphic in $z$ and satisfies the following Riccati equation
\[\chi_z-\frac{\chi^2}2=f.\]
Since $f$ is holomorphic, this Riccati equation has solutions.
So we can find a reciprocal transformation $\Phi$ such that $\Phi(P)$ is local. The $n=2$ case is proved.

The $n=1$ case is similar to the $n=2$ case, we omit the details here.
\end{prf}

\subsection{Locality problem} \label{sec-local}

In this subsection, we investigate the following problem:
{\it Under what conditions is the reciprocal transformation of a Jacobi structure local, i.e. a Hamiltonian structure?}

We fix an invertible $\rho\in\A$, and denote by $\Phi$ the reciprocal transformation w.r.t. $\rho$. The following local functionals are very useful
in this subsection.
\[\Lambda=\int \rho\,dx\in\F,\quad \tilde{\Lambda}=\int 1\,d\tilde{x}\in\tF.\]
It is easy to see that $\Phi(\Lambda)=\tilde{\Lambda}$, $[\tilde{P}, \tilde{\Lambda}]=-\p_{\tilde{\zeta}}(\tilde{P})$.

\begin{lem}\label{local-2}
Let $\tilde{P}\in\tE^p$ with $p\ne2$, then $\tilde{P}$ is local if and only if $\p_{\tilde{\zeta}}(\tilde{P})=0$.
\end{lem}
\begin{prf}
We take a representative of $\tilde{P}=\int \left(\tilde{P}_0+\tilde{\zeta}\,\tilde{X}\right)\,d\tilde{x}$, where $\tilde{P}_0\in\tT_{loc}^p$, 
$\tilde{X}\in\tT_{loc}^{p-1}$. The condition $\p_{\tilde{\zeta}}(\tilde{P})=\int \tilde{X}\,d\tilde{x}=0$ implies that there exists $\tilde{Y}\in\tT^{p-1}$
such that $\tilde{X}=\tp(\tilde{Y})$. Suppose $\tilde{Y}=\tilde{Y}_0+\tilde{Y}_1\tilde{\zeta}$ with $\tilde{Y}_0, \tilde{Y}_1$ local, then the condition
$\tilde{X}\in\tT_{loc}^{p-1}$ implies that $\tilde{Y}_1$ must be a constant. But we have assume $p\ne2$, $\tilde{Y}_1\tilde{\zeta}$ cannot be element of
$\tT^{p-1}$ unless $\tilde{Y}_1=0$, so $\tilde{Y}$ is also local, then
\[\tilde{P}\,\tilde{\sim}\,\tilde{P}_0+\tilde{\zeta}\,\tilde{X}=\tilde{P}_0+\zeta\,\tp(\tilde{Y})\,\tilde{\sim}\,\tilde{P}_0-\tp(\tilde{\zeta})\,\tilde{Y},\]
which is local. The lemma is proved.
\end{prf}

\begin{thm}
Let $P\in\E^p$ with $p\ne2$, then $\Phi(P)$ is local if and only if $[P, \Lambda]=0$.
\end{thm}
\begin{prf}
According to the above lemma, $\Phi(P)$ is local if and only if $\p_{\tilde{\zeta}}\left(\Phi(P)\right)=0$, which is equivalent to $[\Phi(P), \tilde{\Lambda}]=0$,
but we have $[\Phi(P), \tilde{\Lambda}]=\Phi\left([P, \Lambda]\right)$, so the theorem is proved.
\end{prf}

When $p=2$, it is easy to see that the condition $[P, \Lambda]=0$ is still necessary to imply the locality of $\Phi(P)$, but not sufficient.
According to the proof of Lemma \ref{local-2}, we have the following definition.

\begin{dfn}\label{dfn-ch}
Let $P\in\E^2$, and $\rho$ be an invertible differential polynomial satisfying $[P, \Lambda]=0$, then the reciprocal transformation $\Phi(P)$
must take the following form:
\[\Phi(P)=\int \left(\tilde{P}_0+z\,\tilde{u}^{i,1}\tilde{\theta}_i\,\tilde{\zeta}\right)\,d\tilde{x},\]
where $\tilde{P}_0\in\tT^2_{loc}$ and $z\in\R$. The constant $z$ is called the nonlocal charge of the pair $(P, \rho)$, and is denoted by $z(P, \rho)$.
\end{dfn}

\begin{thm} \label{thm-locality}
Let $P\in\E^2$, then $\Phi(P)$ is local if and only if $[P, \Lambda]=0$, $z(P, \rho)=0$.
\end{thm}

In the remaining part of this subsection, we will give a computation formula for the nonlocal charge of a Jacobi structure with a hydrodynamic leading term.

Let $P=(P_0, X)$ be a Jacobi structure, and $\rho$ be an invertible differential polynomial satisfying $[P, \Lambda]=0$, then Lemma \ref{lem-XH} shows
that $[X, \Lambda]=0$, so there exists $\sigma\in\A$ such that
\[\p^s(\delta^i(X))\p_{i,s}(\rho)=\p(\sigma).\]
Suppose the leading term of $P$ is of hydrodynamic type,
\[P_0=\frac12\left(g^{ij}\,\theta_i\,\theta_j'+\Gamma^{ij}_k\,u^{k,1}\,\theta_i\,\theta_j+\cdots\right),\  X=\left(V^i_k\,u^{k,1}+\cdots\right)\theta_i,\]
and suppose the decomposition of $\rho$ and $\sigma$ w.r.t. the gradation \eqref{grad} reads
\[\rho=\rho_0+\cdots,\ \sigma=\sigma_0+\cdots,\]
here $\cdots$ denote the terms with higher degrees, then it is easy to see
\begin{equation}
V^i_k\nabla_i(\rho_0)=\nabla_k(\sigma_0). \label{eq-vrs}
\end{equation}
By straightforward computation, one can obtain
\[[P,\Lambda]=-\int \left(\left(\nabla^i\nabla_k(\rho_0)+\rho_0\,V^i_k+\sigma_0\delta^i_k\right)u^{k,1}+\cdots\right)\theta_i\,dx,\]
if $[P, \Lambda]=0$, then there exists $c\in\R$ such that
\begin{equation}
\nabla^i\nabla_k(\rho_0)+\rho_0\,V^i_k+\sigma_0\delta^i_k=c\,\delta^i_k. \label{eq-nnr}
\end{equation}

On the other hand, Lemma \ref{lem-FP} shows that the nonlocal part of $\Phi(P)$ is
\[(\Phi(P))_{nl}=\int \tilde{\zeta}\,\tilde{V}^i_k\,\tilde{u}^{k,1}\tilde{\theta}_i\,d\tilde{x},\]
where $\tilde{V}^i_k$ reads
\begin{align*}
\tilde{V}^i_k=&\rho_0^2 V^i_k+\rho_0\,\nabla^i\nabla_k\rho_0-\frac12(\nabla^l\rho_0)(\nabla_l\rho_0)\delta^i_k\\
=&\left(\rho_0\left(c-\sigma_0\right)-\frac12(\nabla^l\rho_0)(\nabla_l\rho_0)\right)\delta^i_k,
\end{align*}
so we have
\[(\Phi(P))_{nl}=\int \left(\frac12(\nabla^l\rho_0)(\nabla_l\rho_0)+\rho_0\left(\sigma_0-c\right)\right)
\,\tilde{u}^{i,1}\tilde{\theta}_i\,\tilde{\zeta}\,d\tilde{x},\]
comparing with Definition \ref{dfn-ch}, we obtain
\begin{equation}
z(P, \rho)=\frac12(\nabla^l\rho_0)(\nabla_l\rho_0)+\rho_0\left(\sigma_0-c\right). \label{eq-zpr}
\end{equation}
By using the equations \eqref{eq-vrs} and \eqref{eq-nnr}, one can easily show that the right hand side of \eqref{eq-zpr} is indeed a constant, and is independent of the
choice of $\sigma$, so we get the correct formula for the nonlocal charge of $(P, \rho)$.

\begin{cor} \label{cor-local}
Let $P$ be a quasi-local bivector whose leading term is a  Jacobi structure of hydrodynamic type, suppose $\rho, \rho'$ are two invertible differential polynomials
satisfying $\rho\sim\rho'$, $\Phi, \Phi'$ are the reciprocal transformations w.r.t. $\rho, \rho'$ respectively, then $\Phi(P)$ is local if and only if $\Phi'(P)$
is local.
\end{cor}
\begin{prf}
According to Theorem \ref{thm-locality}, $\Phi(P)$ is local if and only if $[P, \int \rho\,dx]=0$ and $z(P, \rho)=0$.
Since $\rho\sim\rho'$, we have $[P, \int \rho'\,dx]=0$. Note that $z(P, \rho)$ \eqref{eq-zpr} only depends on the degree zero part of $\rho$,
which coincides with the one of $\rho'$, so we also have $z(P, \rho')=0$. Then by using Theorem \ref{thm-locality} again, we obtain that $\Phi'(P)$ is local.
The corollary is proved.
\end{prf}

\subsection{Deformation and cohomology}\label{sec-dfm}

In Section \ref{sec-hydro}, we proved that every Jacobi structure of hydrodynamic type
is equivalent to a Hamiltonian structure of hydrodynamic type modulo reciprocal transformations.
Then by using the classical result of Dubrovin and Novikov  \cite{dn83, dn84} that every  Hamiltonian structure of hydrodynamic type is equivalent to a Hamiltonian structure of
the form 
\begin{equation}
P_{norm}=\int \frac12 \eta^{ij}\theta_i\theta_j^1\,dx
\end{equation}
modulo coordinate transformations, we arrive at the classification of Jacobi
structures of hydrodynamic type: {\em every Jacobi structure of hydrodynamic type is equivalent to a Hamiltonian structure of the form $P_{norm}$.} 

In this subsection, we will generalize the above
classification result to more general Jacobi structures.

Let $P$ be an arbitrary Jacobi structure. Note that $(\E, d_P)$ forms a differential graded Lie algebra (DGLA), where $d_P=\ad_P$ is the adjoint action of $P$:
\[\ad_P:\E\to\E,\ Q\mapsto \ad_P(Q)=[P, Q],\]
so we immediately obtain the following definitions and results along the line of the general philosophy about the relationship between deformation problems
and the cohomologies of the associated DGLAs.

\begin{dfn}\label{dfn-dfm}
Let $P_0\in\E^2$ be a Jacobi structure, we say $P=P_0+Q\in\E^2$ is a deformation of $P_0$ if $P$ is also a Jacobi structure and $\nu(Q)>\nu(P_0)$.
Let $P, P'$ be two deformations of $P_0$, we say $P$ is equivalent to $P'$ if there exists $X\in \E^1$ such that $\nu(X)>0$ and $e^{\ad_X}(P)=P'$.
A deformation $P$ is called trivial if it is equivalent to $P_0$.
\end{dfn}
Note that in the above definition the function $\nu$ is introduced in Definition \ref{dfn-zh}.

\begin{dfn}\label{dfn-idfm}
Let $P\in\E^2$ be a Jacobi structure, $Q\in\E^2$ is called an infinitesimal deformation of $P$ if $\nu(Q)>\nu(P)$ and $d_P(Q)=0$.
Let $Q, Q'$ be two infinitesimal deformations of $P$, they are called equivalent if there exists $X\in\E^1$ such that $\nu(X)>0$ and $Q'-Q=d_P(X)$.
An infinitesimal deformation $Q$ is called trivial if it is equivalent to $0$.
\end{dfn}

\begin{dfn}
Let $P\in\E^2$ be a Jacobi structure, due to  $d_P^2=0$ we have the following complex
\[0 \xrightarrow{\ \ \ } \E^0 \xrightarrow{\ d_P\ } \E^1 \xrightarrow{\ d_P\ } \E^2 \xrightarrow{\ d_P\ } \cdots,\]
its cohomology  is called the Lichnerowicz-Jacobi cohomology of $P$
\begin{equation}
H^i(\E, P)=\frac{\mathrm{Ker}\left(d_P:\E^i\to\E^{i+1}\right)}{\mathrm{Im}\left(d_P:\E^{i-1}\to\E^i\right)}.
\end{equation}
\end{dfn}

If $P$ is homogeneous, i.e. $P=P_{\nu(P)}$, then $H^i(\E, P)$ possesses a gradation
\[H^i(\E, P)=\bigoplus_{d\ge0}H^i_d(\E, P),\ H^i_d(\E, P)=\frac{\mathrm{Ker}\left(d_P:\E^i_d\to\E^{i+1}_{d+\nu(P)}\right)}
{\mathrm{Im}\left(d_P:\E^{i-1}_{d-\nu(P)}\to\E^i_d\right)},\]
where $\E^i_d=\{f \in\E^i|\deg f=d\}$.

\begin{rmk}
Let $\Phi$ be a reciprocal transformation, and $P$ be a Jacobi structure.
Since $\Phi$ preserves the Schouten-Nijenhuis bracket, it induces an isomorphism of complexes $\Phi:(\E, d_P)\to(\tE, d_{\Phi(P)})$, which implies
\[H^*(\E, P)\cong H^*(\tE, \Phi(P)),\]
so reciprocal transformations preserve the Lichnerowicz-Jacobi cohomology. When $P$ is a Jacobi structure of degree zero (see Example \ref{zh-bj-5}),
then $H^*_0(\E, P)$ coincides with the classical Lichnerowicz-Jacobi cohomology defined in \cite{LLMP}. If $\Phi$ is given by a nowhere zero smooth
function on $M$, then $\Phi(P)$ is just a conformal change of $P$ (see Example \ref{zh-bj-5}), so the above isomorphism shows that
the Lichnerowicz-Jacobi cohomology is invariant under conformal changes of the Jacobi structure. This fact was first shown in \cite{LLMP}, our present
proof seems simpler than the one given there.
\end{rmk}

\begin{prp}
Let $P$ be a homogeneous Jacobi structure.

a) The space of equivalence classes of infinitesimal deformations of $P$ is given by $H^2_{>\nu(P)}(\E, P)$.
In particular, every deformation of $P$ is trivial if and only if $H^2_{>\nu(P)}(\E, P)$ vanishes.

b) If $H^3_{>2\nu(P)}(\E, P)$ vanishes, then every infinitesimal deformation can be extended to a genuine deformation,
and the space of equivalence classes of deformations of $P$ is just $H^2_{>\nu(P)}(\E, P)$.
\end{prp}

In the next subsection, we will prove the following theorem.
\begin{thm}\label{thm-trivial}
Let $P$ be a Jacobi structure of hydrodynamic type, then we have
\[H^q_{>0}(\E, P)\cong0, \mbox{ for } q=0, 1, 2, \dots.\]
\end{thm}
From this theorem we immediately have the following corollaries.
\begin{cor}\label{cor-class}
Let $P$ be a Jacobi structure of hydrodynamic type, then every deformation of $P$ can be transformed to $P_{norm}$ by a series of Miura type transformations
and reciprocal transformations.
\end{cor}
\begin{cor}\label{cor-hh}
Let $P$ be a deformation of a Jacobi structure of hydrodynamic type, if an $X\in\E^1_{loc}$ satisfies $[P, X]=0$, then the associated system of evolutionary PDEs has a Jacobi structure, i.e. there exists $H\in\F$ such that
$X=[P, H]$.
\end{cor}

\subsection{Computation of cohomologies}

We fix a Jacobi structure of hydrodynamic type $P\in\E^2$ in this subsection. Note that Miura type transformations and reciprocal transformations preserve the
Lichnerowicz-Jacobi cohomologies, so we can assume $P=\int \frac12 \eta^{ij}\theta_i\theta_j^1\,dx$ without loss of generality.

To every $Q\in\E$, one can associate a first order differential operator $D_Q:\T\to\T$ (see \eqref{zh-09-1}) such that
\[\int D_Q(\alpha)\,dx=[Q, \int \alpha \,dx], \mbox{ for all } \alpha\in\T.\]
For our fixed $P$, the operator $D_P$ reads
\[D_P=\left(\eta^{ij}\theta_i\theta_j^1\right)\p_{\zeta}+\left(\eta^{ij}\theta_j^{s+1}\right)\p_{i,s}\]
and one can prove the following identities
\begin{equation}
D_P^2=0, \quad [D_P, \p]=0.
\end{equation}
Thus we obtain a complex $(\T, D_P)$
\[0 \xrightarrow{\ \ \ } \T^0 \xrightarrow{\ D_P\ } \T^1 \xrightarrow{\ D_P\ } \T^2 \xrightarrow{\ D_P\ } \cdots,\]
and a short exact sequence of complexes
\[0 \xrightarrow{\ \ \ } (\T, D_P) \xrightarrow{\ \p\ } (\T, D_P) \xrightarrow{\ \int\ } (\E, d_P) \xrightarrow{\ \ \ } 0.\]
Then by using the long exact sequence theorem, we obtain the following lemma.
\begin{lem}
For every $d\ge0$, we have the long exact sequence:
\begin{equation}
\cdots \xrightarrow{\ \ \ } H^p_{d-1}(\T) \xrightarrow{\ \p^*\ } H^p_d(\T) \xrightarrow{\ \int^* \ } H^p_d(\E) \xrightarrow{\ \ \ }
H^{p+1}_d(\T)\xrightarrow{\ \ \ } \cdots \label{long1}
\end{equation}
where $H^p_d(\T)=H^p_d(\T, D_P)$, $H^p_d(\E)=H^p_d(\E,d_P)$.
\end{lem}

The complex $(\T, D_P)$ has a subcomplex $(\T_{loc}, D_P')$, where $D_P'$ is the restriction of $D_P$ on $\T_{loc}$:
\[D_P'=\left(\eta^{ij}\theta_j^{s+1}\right)\p_{i,s}.\]
The operator $D_P'$ satisfies $D_P'\p_{\zeta}+\p_{\zeta}D_P=0$, so we have another short exact sequence of complexes
\[0 \xrightarrow{\ \ \ } (\T_{loc}, D_P') \xrightarrow{\ i\ } (\T, D_P) \xrightarrow{\ \p_{\zeta}\ } (\T_{loc}, -D_P') \xrightarrow{\ \ \ } 0,\]
and the following lemma.
\begin{lem}
For every $d\ge0$, we have the long exact sequence:
\begin{equation}
\cdots \xrightarrow{\ \ \ } H^p_d(\T_{loc}) \xrightarrow{\ i^*\ } H^p_d(\T) \xrightarrow{\ \p_{\zeta}^* \ } H^{p-1}_d(\T_{loc}) \xrightarrow{\ \ \ }
H^{p+1}_{d+1}(\T_{loc})\xrightarrow{\ \ \ } \cdots \label{long2}
\end{equation}
where $H^p_d(\T_{loc})=H^p_d(\T_{loc},D_P') \cong H^p_d(\T_{loc},-D_P')$.
\end{lem}

Recall that the algebra $\T_{loc}$ is defined as
\[\T_{loc}=\A\otimes\wedge^*(V_{loc}),\mbox{ where } V_{loc}=\bigoplus_{i,s}\left(\mathbb{R}\theta_i^s\right).\]
We denote $\theta^{i,s}=\eta^{ij}\theta_j^s$, and decompose $V_{loc}$ as $V_{loc}=V_{loc}^{0} \oplus V_{loc}^{>0}$, where
\[V_{loc}^{0}=\bigoplus_{i=1}^n\left(\mathbb{R}\theta^{i,0}\right)\cong\R^n,\ V_{loc}^{>0}=\bigoplus_{i=1}^n\bigoplus_{s>0}\left(\mathbb{R}\theta^{i,s}\right),\]
then $\T_{loc}$ has the following tensor product decomposition
\[\T_{loc}=\T_{loc}^{>0}\otimes\wedge^*(V_{loc}^{0}), \mbox{ where } \T_{loc}^{>0}=\A\otimes\wedge^*(V_{loc}^{>0}),\]
and the differential $D_P'$ becomes $D_P'=\theta^{i,s+1}\p_{i,s}$, which is a differential on $\T_{loc}^{>0}$.
By using the universal coefficient theorem, we obtain the following lemma.
\begin{lem}
\begin{equation}
H^*(\T_{loc}, D_P')\cong H^*(\T_{loc}^{>0}, D_P')\otimes\wedge^*(V_{loc}^{0}). \label{eq-tensor}
\end{equation}
\end{lem}

By regarding $\theta^{i,s+1}$ as $du^{i,s}$, it is easy to see that the complex $(\T_{loc}^{>0}, D_P')$ is just the de Rham complex $(J^{\infty}(M), d_{\mathrm{dR}})$ of
$J^{\infty}(M)$ {\em with coefficients in $\A$}. We introduce a map
\[F:[0,1]\times J^{\infty}(M) \to J^{\infty}(M),\ (t, u^{i,s})\mapsto (t^s\,u^{i,s}),\]
it gives a homotopy from $s_0\circ \pi_{\infty}$ to $\mathrm{Id}_{J^{\infty}(M)}$, where $s_0$ and $\pi_{\infty}$ is the zero section and the
structure projection of the bundle $J^{\infty}(M)$ respectively, so the induced maps
\[\left(s_0\circ \pi_{\infty}\right)^*, \left(\mathrm{Id}_{J^{\infty}(M)}\right)^*: H^*(J^{\infty}(M), d_{\mathrm{dR}})\to H^*(J^{\infty}(M), d_{\mathrm{dR}})\]
must be same. But note that when $d>0$, the induced map
\[\left(s_0\circ \pi_{\infty}\right)^*: H^*_d(J^{\infty}(M), d_{\mathrm{dR}})\to H^*_d(J^{\infty}(M), d_{\mathrm{dR}})\]
is zero, so we obtain
\begin{equation}
H^*(\T_{loc}^{>0}, D_P')\cong H^*_0(J^{\infty}(M), d_{\mathrm{dR}}) \cong H^*_{\mathrm{dR}}(M)\cong \R, \label{eq-derham}
\end{equation}
since we have assumed that $M$ is contractible.

\vskip 1em

\begin{prfn}{Theorem \ref{thm-trivial}}
From Equations \eqref{eq-derham} and \eqref{eq-tensor} we obtain
\[H^p_d(\T_{loc}, D_P')\cong\left\{\begin{array}{cc} \wedge^p(\R^n), & d=0, \\ 0, & d>0,\end{array}\right.\]
so the long exact sequence \eqref{long2} implies $H^*_{>0}(\T)\cong0$, then the theorem is proved by applying the long exact sequence \eqref{long1}.
\end{prfn}

\subsection{Bi-Jacobi structures}\label{sec-35}

In this subsection we introduce the notion of compatibility of two Jacobi structures, and use the Lenard-Magri recursion scheme to study integrability of a system which possesses a pair of compatible Jacobi structures, i.e., a bi-Jacobi structure.
We also consider properties of a bi-Jacobi structure under reciprocal transformations.

\begin{dfn}We say that two Jacobi structures
$P, Q\in\E^2$ are compatible if  $[P, Q]=0$, and call a compatible pair of Jacobi
structures a bi-Jacobi structure. We say that an evolutionary PDE $X\in\E^1_{loc}$
possesses a bi-Jacobi structure if there exists a bi-Jacobi structure $(P, Q)$ and two local functionals $H_0, H_1$ such that
\begin{equation}\label{zh-09-21}
X=[P, H_1]=[Q,H_0].
\end{equation}
\end{dfn}

In what follows we assume that both leading terms of $P$ and $Q$ are of hydrodynamic type.
According to Corollary \ref{cor-class}, there exists a reciprocal transformation $\Phi$ such that $\Phi(P)$ is local, then $\Phi(X)$ is also local, this is due to the following identity:
\[\Phi(X)=\Phi\left([P,H_1]\right)=[\Phi(P), \Phi(H_1)].\]
So without loss of generality we assume that $P$ is local.

\begin{lem}
For a given evolutionary PDE \eqref{zh-09-21} that possesses a bi-Jacobi structure,
there exists local functionals $H_2, H_3, \dots \in\F$ such that
\begin{equation}\label{zh-09-21-2}
[P, H_{k+1}]=[Q, H_k],\ k=1, 2, \dots
\end{equation}
and moreover, for any $k, l\ge0$ we have
\begin{equation}
P(H_k, H_l)=0,\ Q(H_k, H_l)=0. \label{eq-lmrs}
\end{equation}
\end{lem}
\begin{prf}
By using the Jacobi identity and the second equality of \eqref{zh-09-21} we have
\[ [P,[Q,H_1]]=-[Q,[P,H_1]]-[H_1,[P,Q]]=-[Q,[P,H_1]]=-[Q,[Q,H_0]]=0,\]
so from Corollary \ref{cor-hh} it follows the existence of a local functional $H_2$ such that
\[[P,H_2]=[Q,H_1].\]
By repeating this procedure we can obtain local functionals $H_2, H_3,\dots$ satisfying
\eqref{zh-09-21-2}.

The equations \eqref{eq-lmrs} can be proved by using the recursion relation \eqref{zh-09-21-2} and Magri's trick given in \cite{magri}.
\end{prf}

We denote $X_k=[P, H_k]$, then $X_k$'s are all local, and for any $k, l\ge0$ we have
\[[X_k, H_l]=0,\ [X_k, X_l]=0.\]
So in the case when $\{H_k\}_{k=0}^{\infty}$ are linearly independent we obtain an integrable hierarchy $\{X_k\}_{k=1}^{\infty}$ which has infinitely many common
conserved quantities $\{H_k\}_{k=0}^{\infty}$.

\vskip 1em

The notion of bi-Jacobi structure is a generalization of the one of bihamiltonian structure, the above lemma shows its close relation with integrable
systems. An important problem concerning such a relation is the problem of classification of  the orbit space of bi-Jacobi structures
under the action of Miura type transformations and reciprocal transformations. In the local case, i.e. when both $P$ and $Q$ are local, the classification problem
has been partially solved in \cite{LZ-1, DLZ-1}.  The main results of these papers is the 
introduction of the concept of central invariants for any deformation of a semisimple bihamiltonian structures of hydrodynamic
type, and proof of the fact that two deformations are equivalent if and only if they possess same central invariants. In the general nonlocal case the classification
problem is far from being solved.

In what follows of this subsection, we proceed to consider properties of 
a bi-Jacobi structure under reciprocal transformations. To this end we first recall some definitions.

\begin{dfn}
Let $(P_0, Q_0)$ be a bihamiltonian structure of hydrodynamic type
\[P_0\sim\frac12\left(g_1^{ij}\theta_i\theta_j^1+\Gamma^{ij}_{1,k}u^{k,1}\theta_i\theta_j\right),\
Q_0\sim\frac12\left(g_2^{ij}\theta_i\theta_j^1+\Gamma^{ij}_{2,k}u^{k,1}\theta_i\theta_j\right).\]
We say that $(P_0, Q_0)$ is semisimple if the roots of the characteristic equation
\begin{equation}
\det\left(g_2^{ij}-\lambda\,g_1^{ij}\right)=0 \label{eq-chara}
\end{equation}
are distinct and not constant.
\end{dfn}

According to the result of Ferapontov \cite{fera-2}, the property of semisimplicity implies that the roots of the characteristic equation \eqref{eq-chara}
can be used as local coordinates near every point of $M$, we denote them by $\lambda_1, \dots, \lambda_n$ and call them the canonical coordinates of the bihamiltonian
structure $(P_0, Q_0)$ \cite{LZ-1, DLZ-1}. The two metrics become diagonal under the canonical coordinates,
\[g_1^{ij}(\lambda)=\delta^{ij}\, f_i(\lambda),\ g_2^{ij}(\lambda)=\delta^{ij}\, \lambda_i\, f_i(\lambda).\]
From now on in this subsection, we will assume that  $(P_0, Q_0)$ is a semisimple bihamiltonian structure represented in its canonical coordinates.

Let $(P, Q)$ be a deformation of $(P_0, Q_0)$, we introduce the following tensors:
\begin{equation}
A^{ij}_m=\left(\p^j_m\delta^i(P)\right)_0,\ B^{ij}_m=\left(\p^j_m\delta^i(Q)\right)_0,\ m=0, 1, 2, \dots \label{eq-AB}
\end{equation}
where $(W)_0$ denotes the degree zero part of $W\in\A$. For example, $A^{ij}_0=B^{ij}_0=0$, $A^{ij}_1=g^{ij}_1$, $B^{ij}_1=g^{ij}_2$, and so on.

\begin{dfn}
The central invariants of $(P, Q)$ are defined by
\begin{equation}\label{dfn-ci}
c_i(\lambda)=\frac1{3\,f_i^2} \left(B^{ii}_3- \lambda_i A^{ii}_3+\sum_{k\ne i}\frac{(B^{ki}_2-\lambda_i A^{ki}_2)^2}{f_k
(\lambda_k-\lambda_i)}\right),\ i=1, \dots, n.
\end{equation}
It can be shown that $c_i$ only depends on $\lambda_i$.
\end{dfn}

We have the following important property of the central invariants:
\begin{thm}[\cite{LZ-1, DLZ-1}]\label{thm-civ0}
Let $(P, Q),\ (P', Q')$ be two deformations of $(P_0, Q_0)$, and $c_i, c_i'\ (i=1, \dots, n)$
be their central invariants respectively, then $(P, Q)$ is equivalent to $(P', Q')$ via a Miura type transformation if and only if $c_i=c_i'\ (i=1, \dots, n)$.
\end{thm}

The main result of this subsection is the following theorem.
\begin{thm} \label{thm-civ}
Let the bihamiltonian structure $(P, Q)$ be a deformation of $(P_0,Q_0)$,
and $\Phi$ be a reciprocal transformation such that both
$\Phi(P)$ and $\Phi(Q)$ are local. Then $(\Phi(P), \Phi(Q))$ is also a deformation of a semisimple bihamiltonian structure of hydrodynamic
type which has the same system of canonical coordinates as $(P_0,Q_0)$, and its central invariants coincide with that of $(P, Q)$.
\end{thm}
\begin{prf}
Since reciprocal transformations induce conformal changes of the metrics $g_1$ and $g_2$, it is easy to see that the canonical coordinates are preserved under
such transformations. We denote by $\tilde{A}^{ij}_m$, $\tilde{B}^{ij}_m$ the tensors \eqref{eq-AB} associated to $(\Phi(P), \Phi(Q))$,
then by using the formula \eqref{trsp} we can easily obtain
\[\tilde{A}^{ij}_m=\rho^{m+1}A^{ij}_m,\ \tilde{B}^{ij}_m=\rho^{m+1}B^{ij}_m.\]
By using these relations the theorem follows from the definition of the central invariants immediately.
\end{prf}

\begin{rmk}
In Theorem \ref{thm-civ} the condition that $\Phi(P)$, $\Phi(Q)$ are local is very important. There are some important bihamiltonian integrable hierarchies
which are equivalent via reciprocal transformations but have different central invariants. They are not counter-examples to Theorem \ref{thm-civ}
because these integrable hierarchies possess more than two Jacobi structures, after the reciprocal transformation some of their local Jacobi
structures may become nonlocal, and some nonlocal ones may become local. So the resulting hierarchy may still possess a (local) bihamiltonian structure, but this
bihamiltonian structure is not the one that is transformed from the original local one, so the theorem cannot apply to such cases.
\end{rmk}

\begin{cor}\label{cor-sim}
Let $(P, Q)$ be a bihamiltonian structure whose leading term is semi\-simple and of hydrodynamic type. For two invertible differential polynomials $\rho, \rho'\in\A$ satisfying $\rho\sim\rho'$, let $\Phi, \Phi'$ be the reciprocal transformations associated to them respectively. Then if both $\Phi(P), \Phi(Q)$ are local,
there exists a Miura type transformation $g$ such that
\[g\left(\Phi(P)\right)=\Phi'(P),\ g\left(\Phi(Q)\right)=\Phi'(Q).\]
\end{cor}
\begin{prf}
According to Corollary \ref{cor-local}, $\Phi'(P)$ and $\Phi'(Q)$ are both local.
Suppose $\rho'=\rho+\p F$, then $\eta=\rho'/\rho=1+\tp\,F$. Let $\Psi$ be the reciprocal transformation w.r.t. $\eta$, then we have
\[\Psi\left(\Phi(P)\right)=\Phi'(P),\ \Psi\left(\Phi(Q)\right)=\Phi'(Q),\]
so $\Psi$ is a reciprocal transformation converting a local bihamiltonian structure to a local bihamiltonian structure, so the central invariants of
$\left(\Phi(P), \Phi(Q)\right)$ and $\left(\Phi(P'), \Phi(Q')\right)$ must coincide. On the other hand, it is obvious that the bihamiltonian structures
$\left(\Phi(P), \Phi(Q)\right)$ and $\left(\Phi(P'), \Phi(Q')\right)$ possess same leading terms, so they are deformations of a same semisimple bihamiltonian
structure of hydrodynamic type. Thus according to Theorem \ref{thm-civ0} there exists a Miura type transformation $g$ such that
\[g\left(\Phi(P)\right)=\Phi'(P),\ g\left(\Phi(Q)\right)=\Phi'(Q),\]
the corollary is proved.
\end{prf}

We note that the above theorem  only ensures
the existence of $g$, we do not have the explicit form of $g$ in general.

\section{Examples}

\subsection{The KdV equation v.s. the Camassa-Holm equation}

In this subsection, we investigate the relation between the Jacobi structures of the KdV equation and the Camassa-Holm equation.

The KdV equation
\begin{equation}
u_t=X=6\,u\,u_x-u_{xxx}, \label{eq-kdv}
\end{equation}
possesses two compatible local Hamiltonian structures with the following Hamiltonian operators \cite{gardner, zf, magri}:
\[P_0=\p,\ P_1=u\,\p+\frac12\,u_x-\frac14\,\p^3.\]

Denote $R=P_1 \circ P_0^{-1}$, and $P_k=R^k\circ P_0$, then $R$ is a hereditary strong symmetry \cite{FF}, and all $P_k$'s are Hamiltonian operators
which are, however, 
nonlocal for $k\ge 2$.
Among these nonlocal Hamiltonian operators, there are two operators that correspond to Jacobi structures:
\begin{align*}
P_2=&u^2\,\p+u\,u_x-\frac14\,u_x\,\p^{-1}\,u_x\\
&-\left(\frac12\,u\,\p^3+\frac34\,u_x\,\p^2+\frac12\,u_{xx}\,\p+\frac18\,u_{xxx}\right)+\frac1{16}\,\p^5, \\
P_3=&u^3\,\p+\frac32\,u^2\,u_x\\&-\frac1{16}\,\left(6u\,u_x-u_{xxx}\right)\p^{-1}\,u_x-\frac1{16}\,u_x\,\p^{-1}\left(6u\,u_x-u_{xxx}\right)\\
&-\left(\frac34\,u^2\,\p^3+\frac94\,u\,u_x\,\p^2+\left(\frac32\,u\,u_{xx}+\frac{17}{16}\,\,u_x^2\right)\p\right.\\
&\left.+\left(\frac38\,u\,u_{xx}+\frac{11}{16}\,u_x\,u_{xx}\right)\right)+\left(\frac3{16}\,u\,\p^5+\frac{15}{32}\,u_x\,\p^4\right.\\
&\left.+\frac58\,u_{xx}\p^3+\frac{15}{32}\,u_{xxx}\p^2+\frac{3}{16}\,u_{4x}\p+\frac1{32}\,u_{5x}\right)-\frac1{64}\p^7
\end{align*}
We denote by $Q_i\ (i=0,1,2,3)$ the quasi-local bivectors corresponding to $P_i\ (i=0,1,2,3)$
\[Q_i=\frac12\,\int\theta\,P_i(\theta)\,dx,\]
then they read
\begin{align*}
Q_0=&\frac12\,\int \theta\theta'\,dx, \\
Q_1=&\frac12\,\int \left(u\,\theta\theta'+\frac32\theta\theta'''\right)\,dx,\\
Q_2=&\frac12\,\int \left(u^2\,\theta\theta'+\frac14\,u_x\,\theta\,\zeta+\cdots\right)\,dx,\\
Q_3=&\frac12\,\int \left(u^3\,\theta\theta'+\frac34\,u\,u_x\,\theta\,\zeta+\cdots\right)\,dx,
\end{align*}
where we omitted some terms with higher degrees.
Here an interesting problem arises: does there exist more Jacobi structures other than $Q_0, Q_1, Q_2, Q_3$ for the KdV equation \eqref{eq-kdv}?
We conjecture that the answer is no, and we will investigate this problem in a separate publication \cite{LZ-3}.

Let $\rho\in\A$ be an invertible differential polynomial, and $\Phi$ be the reciprocal transformation w.r.t. $\rho$.
It is easy to see that to ensure 
the locality of $\Phi(X)$ the differential polynomial  $\rho$ must be a
conserved density of the KdV hierarchy.
From the quasi-triviality of the KdV hierarchy \cite{LZ-2} it follows that the conserved densities are uniquely determined by their leading terms
(up to a total differential). We assume 
that
the leading term of $\rho$ is $u^{\kappa}\ (\kappa\in\R)$, then by using Lemma \ref{lem-FP} and the 
results of Section \ref{sec-local} one can obtain the following assertions:
\begin{itemize}
\item[(a)] $\Phi(Q_0)$ is local $\Leftrightarrow$ $\kappa=0,\ 2$; \quad (b) $\Phi(Q_1)$ is local $\Leftrightarrow$ $\kappa=0,\ 1$;
\item[(c)] $\Phi(Q_2)$ is local $\Leftrightarrow$ $\kappa=\frac12,\ -\frac12$; \quad (d) $\Phi(Q_3)$ is local $\Leftrightarrow$ $\kappa=\frac12,\ -\frac32$.
\end{itemize}

The case $\kappa=0$ is trivial. 
When $\kappa=2$, we can choose $\rho=u^2$, the transformed equation $\Phi(X)$ is analyzed in \cite{LZ-2} and \cite{LWZ-1}.
It is shown in \cite{LWZ-1} that $\Phi(X)$ has only one local Hamiltonian structure. Now we know that it also has (at least) three Jacobi structures.
The case when $\kappa=1$ is similar.

We are interested in the case 
when
$\kappa=\frac12$. 
For this we first need to find out the full expression of $\rho$.
It is well known
that the following Miura transformation
\[u=v^2+\,v_x\]
converts the KdV
equation \eqref{eq-kdv} to the modified KdV equation
\[v_t=6\,v^2\,v_x-v_{xxx}=\left(2\,v^3-v_{xx}\right)_x,\]
so we can choose the conserved density $\rho=2\,v$. Here the factor $2$ is not essential, we add it just for the convenience of relating the KdV equation with the standard Camassa-Holm equation.
After fixing this $\rho$,
one obtains $\Phi(Q_i)\ (i=0, 1, 2, 3)$ with $\Phi(Q_2)$ and $\Phi(Q_3)$ being local.

\begin{thm}
There exists a Miura type transformation $g$ such that
\[\left(g\left(\Phi(Q_2)\right),g\left(\Phi(Q_3)\right)\right)\]
gives the bihamiltonian structure of 
the
Camassa-Holm hierarchy \cite{CH, CHH, Fu, FF}.
\end{thm}
\begin{prf}
We recall some facts about the Camassa-Holm (CH) equation first. The CH equation
\begin{equation}
m_s=2\,m\,w_y+m_y\,w,\quad m=w-w_{yy} \label{eq-ch}
\end{equation}
has two local Hamiltonian structures
\[J_0=\p_y-\p_y^3,\quad J_1=m\,\p_y+\frac12\,m_y\]
and infinitely many nonlocal Hamiltonian structures $J_k=\left(J_1J_0^{-1}\right)^k J_0$
among  which
$J_2$ and $J_3$ correspond to Jacobi structures.
The CH hierarchy can be obtained by using the recursion operator $J_1 J_0^{-1}$.

The simplest way to describe the CH hierarchy is to consider the isospectral problem
\begin{equation}
\phi_{yy}=\left(\frac14-m\,\lambda\right)\phi,\quad \phi_s=-B(\lambda)\phi_y+\frac12B_y(\lambda)\phi, \label{eq-isosp}
\end{equation}
where $B(\lambda)$ is a Laurent polynomial of the spectral parameter $\lambda$. For example, when
\[B(\lambda)=w+\frac1{2\lambda},\]
we obtain
the CH equation \eqref{eq-ch} from the 
compatibility
condition $\phi_{yys}=\phi_{syy}$.

The CH hierarchy has a conserved density $\tilde{\rho}=\sqrt{m}$, 
we denote the associated reciprocal transformation by $\tilde{\Phi}$.
Let
$\psi=m^{\frac14}\,\phi$, then the $y$-part of the isospectral problem \eqref{eq-isosp} becomes
\begin{equation}
\psi_{xx}=\left(u-\lambda\right)\psi, \label{eq-lax}
\end{equation}
where $x=\tilde{y}$ is defined by 
\[d\tilde{y}=\tilde{\rho} \,dy+w \sqrt{m}\, ds,\quad d\tilde{s}=ds\]
and the function $u$ reads
\begin{equation}
u=\frac1{4\,m}+\frac{m_{xx}}{4\,m}-\frac{3\,m_x^2}{16\,m^2}. \label{eq-chkdv}
\end{equation}

Note that the equation \eqref{eq-lax} is exactly the $x$-part of the isospectral problem of the KdV equation \eqref{eq-kdv}.
In fact, if we take the $t$-part as
\[\psi_t=\left(4\,\lambda+2\,u\right)\psi_x-u_x\,\psi,\]
then the compatibility
condition $\psi_{xxt}=\psi_{txx}$ gives us the KdV equation \eqref{eq-kdv}.
So we see that the reciprocal transformation $\tilde{\Phi}$ together with the Miura type transformation \eqref{eq-chkdv}
transforms the CH hierarchy to the negative flows of the KdV hierarchy.
Let $\Phi'$ be the inverse of $\tilde{\Phi}$, then the reciprocal transformation $\Phi'$ is defined by
\[\rho'=\tilde{\rho}^{-1}=2\sqrt{u}+\cdots.\]
Note that $\rho'$ must be 
a conserved
density of the KdV hierarchy, so by using the quasi-triviality property of the KdV hierarchy \cite{LZ-2} we have $\rho'\sim \rho$.

We denote the quasi-local bivectors corresponding to $J_i\ (i=0, 1, 2, 3)$ by $K_i\ (i=0, 1, 2, 3)$.
Since $\rho'$ is a conserved density with leading term $2\sqrt{u}$, we know from Corollary \ref{cor-local} that
the Jacobi structures $\Phi'(Q_2), \Phi'(Q_3)$
are also local. 
By a
straightforward computation, one can show that the leading terms of $\Phi'(Q_2), \Phi'(Q_3)$ coincide with 
that
of $K_1, K_0$
(after the Miura type transformation \eqref{eq-chkdv}). On the other hand, the main results of \cite{LWZ-1} 
imply that the Hamiltonian structures
of the CH equation are uniquely determined by their leading terms, so we have
\[\Phi'(Q_2)=K_1,\quad  \Phi'(Q_3)=K_0.\]
A similar argument shows that
\[\tilde{\Phi}(K_2)=Q_1,\quad  \tilde{\Phi}(K_3)=Q_0,\]
thus we have $\Phi'(Q_0)=K_3, \Phi'(Q_1)=K_2$.

Finally, due to Corollary \ref{cor-sim} there exists a Miura type transformation $g$ such that
\[g(\Phi(Q_2))=\Phi'(Q_2),\quad  g(\Phi(Q_3))=\Phi'(Q_3),\]
the theorem is proved.
\end{prf}

The relation between the Jacobi structures of the KdV equation 
and that of
the CH equation has a multi-component generalization, which
relates the Jacobi structures of the $r$-KdV-CH hierarchy associated to the parameters 
$(a_0, \dots, a_r) \in\left(\R^{r+1}\right)^{\times}$
with the Jacobi structures of the $r$-KdV-CH hierarchy associated to the parameters
$(a_r, \dots, a_0)$, see \cite{CLZ-2} for more details.

\subsection{Bihamiltonian structures associated to Frobenius manifolds}
In this subsection we consider reciprocal transformation of the bihamiltonian structure defined on the jet space
of a Frobenius manifold. The notion of Frobenius manifold was introduced by Dubrovin\cite{D-2, D-1} as a geometrical interpretation of the
WDVV equations that arose in the study of 2D topological field theory \cite{DVV, Witten1}. On any Frobenius manifold $M$ there is defined a flat metric
$\eta$, and on each of its tangent spaces there is a Frobenius algebra structure which depends smoothly on the point of the manifold.
The structure constants of the Frobenius algebras can be represented by a function $F$, called the potential of the
Frobenius manifold, of the flat coordinates $v^1,\dots, v^n$ of $M$. These flat coordinates can be chosen in such a way that
the vector field $e_1=\frac{\p}{\p v^1}$ gives the unity element of the Frobenius algebras. Denote by $\eta_{ij}$ the components of the flat metric in the flat coordinates,
and $(\eta^{ij})=(\eta_{ij})^{-1}$, then the operation of multiplication of the Frobenius algebras on the tangent spaces is given by
\[ \frac{\p}{\p v^i} \cdot \frac{\p}{\p v^j}=c^k_{ij} \frac{\p}{\p v^k},\quad
c_{ij}^k=\eta^{kl}\frac{\p^3 F(v)}{\p v^l\p v^i\p v^j}.\]
Note that we have $c^k_{1j}=\delta^k_j$. The associativity condition of the Frobenius algebras yields a system of nonlinear PDEs for the potential $F$. This function
also satisfies certain quasi-homogeneity condition which is represented by 
an Euler vector field $E$ of the form
\begin{equation}\label{zh-10-7-2}
E=\sum_{i=1}^n E^i\frac{\p}{\p v^i}=\sum_{i=1}^n \left(d_i v^i+r_i\right)\frac{\p}{\p v^i},\quad {\rm{with}}\ d_1=1
\end{equation}
satisfying
\begin{equation}\label{zh-10-7}
\p_E F=(3-d) F+\frac12 A_{ij}\, v^i v^j+ B_i v^i+C.
\end{equation}
Here $d_i, r_i, d, A_{ij}, B_i, C$ are some constants, $A_{1i}=A_{i1}=
\sum_{l=1}^n \eta_{il}\, r_l$, and the constant $d$ is called the charge of the Frobenius manifold. The Euler vector field is so normalized that $r_k=0$ when $d_k\ne 0$.

In what follows of this subsection we assume that the flat coordinates can be normalized so that the metric $\eta$ takes the form
\begin{equation}\label{zh-09-24-1}
\eta_{ij}=\delta_{i,n+1-j},
\end{equation}
and the potential $F(v)$ has the expression
\[F(v)=\frac12 (v^1)^2 v^n+\frac12\,\sum_{k,l=2}^{n-1} \eta_{kl}\, v^1 v^k v^l+f(v^2,\dots, v^n).\]
Then the constants $d_i$ satisfy the following duality property:
\[d_i+d_{n+1-i}=2-d.\]
In the cases when $d=1$ and $d=2$, we also impose the additional conditions $r_n=0$ in \eqref{zh-10-7-2} and $B_1=0$ in \eqref{zh-10-7} respectively.

The flat metric $\eta$ yields a Hamiltonian structure of hydrodynamic type on the jet space of $M$ given by the
Hamiltonian operator
\begin{equation}
P_1^{ij}=\eta^{ij} \p_x,\quad i,j=1,\dots,n.
\end{equation}
There is another compatible Hamiltonian structure which is associated to the intersection form $(g^{ij})$ of the
Frobenius manifold  defined by
\begin{equation}
g^{ij}=E^k c^{ij}_k,\quad c^{ij}_k=\eta^{il} c^j_{lk}.
\end{equation}
The Hamiltonian operator has the form
\begin{equation}
P_2^{ij}=g^{ij}(v)\p_x+\Gamma^{ij}_k(v)\,v^k_x,
\end{equation}
where $\Gamma^{ij}_k=-g^{il}\Gamma^j_{lk}$ are the contravariant components of the Levi-Civita connection of the metric $(g_{ij})=(g^{ij})^{-1}$ defined on the open subset of $M$ where $ \det(g^{ij}(v))\ne 0$.

Let us consider the reciprocal transformation $\Phi$ w.r.t. the function $\rho=v^n$ of the bihamiltonian structure $(P_1, P_2)$.
By using Lemma \ref{lem-FP} we know that the bihamiltonian structure is transformed to a bi-Jacobi structure in general.
However in this particular case the the bihamiltonian structure is transformed to a local bi-Jacobi structure, i.e. again a bihamiltonian structure.
In the coordinates $v^1,\dots, v^n$ this bihamiltonian structure is given by
\begin{equation}
\Phi(P_1)^{ij}=\tilde{\eta}^{ij}\p_{\tilde{x}}+\tilde{\Gamma}^{ij}_{1,k}\, v^k_{\tilde{x}},\quad
\Phi(P_2)^{ij}=\tilde{g}^{ij}\p_{\tilde{x}}+\tilde{\Gamma}^{ij}_{k,2}\, v^k_{\tilde{x}},
\end{equation}
where $\tilde{\eta}^{ij}=\rho^2 \eta^{ij},\ \tilde{g}^{ij}=\rho^2 g^{ij}$, and $\tilde{\Gamma}^{ij}_{1,k}, \tilde{\Gamma}^{ij}_{2,k}$ are the contravariant
components of the Levi-Civita connection of  $(\tilde{\eta}_{ij}),\ (\tilde{g}_{ij})$ respectively.

Note that the metrics $\tilde{\eta}$ and $\tilde{g}$ are also flat. It is easy to see
that the first metric $\tilde{\eta}$ has the following system of flat coordinates:
\begin{equation}\label{zh-09-25-1}
\tilde{v}^1=\frac12\,\frac{\eta_{kl} v^k v^l}{v^n},\quad
\tilde{v}^i=\frac{v^i}{v^n} \ \ (i\ne 1,n),\quad \tilde{v}^n=-\frac1{v^n}.
\end{equation}
In these flat coordinates the metric $\tilde{\eta}$ has the contravariant components
\begin{equation}\label{zh-09-25-2}
\tilde{\eta}^{ij}(\tilde{v})=\eta^{ij}.
\end{equation}

It is noticed in \cite{XZ} that the metrics $\tilde{\eta}(\tilde{v})$ and $\tilde{g}(\tilde{v})$ coincide with the flat metric and
the intersection form of the Frobenius manifold $\tilde{M}$ with the potential
\begin{equation}\label{zh-09-25-3}
\tilde{F}(\tilde{v})=\left(\tilde{v}^n\right)^2 F(v)+\frac12\, \eta_{kl}\, \tilde{v}^1 \tilde{v}^k \tilde{v}^l.
\end{equation}
Here in the r.h.s. of the above formula the variables $v^1,\dots, v^n$ of $F(v)$ are understood to be dependent on $\tilde{v}^1,\dots, \tilde{v}^n$ via \eqref{zh-09-25-1}.
The Euler vector field $\tilde{E}=\sum_{i=1}^n \tilde{E}^i \frac{\p}{\p \tilde{v}^i}$ of the Frobenius manifold $\tilde{M}$ coincides with $E$, i.e. 
\[\tilde{E}=\sum_{i=1}^n \tilde{E}^i \frac{\p}{\p \tilde{v}^i}=
\sum_{i=1}^n E^i \frac{\p}{\p v^i}.\]
and the charge of $\tilde{M}$ is $\tilde{d}=2-d$.

The construction of the Frobenius manifold $\tilde{M}$ from $M$ by the above formula \eqref{zh-09-25-1}--\eqref{zh-09-25-3}
was given by Dubrovin in \cite{D-1}, it corresponds to a special symmetry of the
WDVV equations called the inversion symmetry. The above fact gives a natural interpretation of this symmetry of the WDVV equations in terms of a special reciprocal transformation of the associated bihamltonian structure on the jet space of the Frobenius manifold. Note that under certain additional conditions a bihamiltonian structure of hydrodynamic type determines a Frobenius manifold \cite{du-taniguchi}. In \cite{XZ}
the transformation rule of the associated bihamiltonian integrable hierarchy, called the principal hierarchy of the Frobenius manifold, under the above reciprocal transformation is given, it shows that the reciprocal transformation induces a Legendre transformation of the tau function of the principal hierarchy.

\section{Conclusion}
In the preceding sections we introduce the notion of infinite dimensional Jacobi structures which provides a geometrical characterization of an important class of
nonlocal Hamiltonian structures. The invariance of the space of Jacobi structures under reciprocal transformations enables one to study canonical forms of Jacobi
structures under Miura type transformations and reciprocal transformations, and to discover possible hidden relations between some important evolutionary PDEs via
the study of their local and nonlocal Hamiltonian structures. For a Jacobi structure with leading term of hydrodynamic type, we obtain its canonical form by
computing the associated Lichnerowicz-Jacobi cohomologies.

From the point view of integrable evolutionary PDEs, our next important task is to consider the problem of classification of bi-Jacobi structures with
leading terms of hydrodynamic type. As we have shown in Section \ref{sec-35}, under certain appropriate conditions an evolutionary PDE which has a bi-Jacobi
structure possesses an infinite number of commuting flows and conservation laws. Thus bi-Jacobi structures just like bihamiltonian structures are closely related to
integrable evolutionary PDEs. To classify bi-Jacobi structures of hydrodynamic type, one must classify first their leading terms, i.e. pairs of compatible conformally
flat metrics, which is a pure differential geometric problem; then the classification problem becomes a deformation problem, which is controlled by
some {\it bi-Jacobi cohomologies}. The computation of these bi-Jacobi cohomologies is highly nontrivial, new techniques are expected.

The reciprocal transformations that we consider in this paper only change the spatial variable $x$ of a system of evolutionary PDEs
$u^i_t=K^i(u,u_x,\dots),\ i=1,\dots,n$. It is also important to consider more general class of reciprocal transformations of the form
\begin{equation}\label{zh-10-7-3}
d \tilde{x}=\rho_1 dx+\sigma_1 dt,\quad d\tilde{t}= \rho_2 dx+\sigma_2 dt,
\end{equation}
they change both the time and spatial variables. Here $\frac{\p \rho_i}{\p t}=\frac{\p\sigma_i}{\p x},\ i=1,2$ are two conservation laws of the evolutionary PDEs
which satisfy certain non-degeneracy conditions. Such reciprocal transformations can be represented as the composition of two transformations of the form
\eqref{eq-reci} and one of the form
\begin{equation}\label{zh-10-04}
d\tilde{x}=dt,\quad d\tilde{t}=d x,
\end{equation}
i.e. a transformation that exchanges the role of the time and spatial variables of the evolutionary PDEs. The transformation rule of the Hamiltonian structures of
a system of Hamiltonian evolutionary PDEs under the reciprocal transformation \eqref{zh-10-04} was considered in \cite{pavlov}, where it was shown that for
Hamiltonian systems of hydrodynamic type the reciprocal transformation \eqref{zh-10-04} preserves the forms of the systems and their Hamiltonian structures. It was
further shown in \cite{XTZ} that the reciprocal transformation \eqref{zh-10-04} also preserves the bihamiltonian structures of bihamiltonian systems of hydrodynamic
type. We conjecture that for a general Hamiltonian (resp. bihamiltonian) system with dispersion terms the reciprocal transformation \eqref{zh-10-04} preserves the
localilty property of the Hamiltonian and bihamiltonian structures. If this conjecture holds true, we could consider the problem of classification of bihamiltonian
structures and the associated hierarchies of integrable evolutionary PDEs under Miura type transformations and also reciprocal transformations which preserve the locality
property. It is expected that the moduli spaces that arise in such a classification scheme are much simple than that of the bihamiltonian integrable
hierarchies under the action of Miura type transformations. A preliminary study suggests that every hierarchy of scalar evolutionary PDEs with leading terms of
hydrodynamic type and a bihamiltonian structure is equivalent to the usual KdV hierarchy under Miura type transformations and reciprocal transformations
of the form \eqref{zh-10-7-3}.  We will discuss this issue in separate publications.

The action of reciprocal transformation \eqref{zh-10-04} on general Jacobi structures is still an open problem. The results of \cite{FP} suggest that
reciprocal transformations of the form \eqref{zh-10-04} may convert Jacobi structures to more general nonlocal Hamiltonian structures which are not
rigorously defined yet. In a recent paper \cite{KKVV}, Kersten {\it et al} sketch out a new geometric framework to construct Hamiltonian operators for generic,
non-evolutionary PDEs, we hope that their approach will be helpful to the study of this problem.

The notion of Jacobi structures defined in the present paper is only for $1+1$ dimensional evolutionary PDEs. It is not hard to generalize it to higher spatial dimensional
cases. To do so one only needs to replace the jet space $J^\infty(M)$ used in Section \ref{sec-2} by more general jet spaces. Some results of the present paper still hold
true for Jacobi structures in higher spatial dimensional cases, while some results do not. There are many interesting open problems in this aspect.

\vskip0.6 truecm
\noindent{\bf Acknowledgments.}
The authors thank Boris Dubrovin for his interests in and advises on the work of the paper.
Part of the work was done in SISSA during the authors' visit to it, they thank its hospitality and financial support.
The work is partially supported by the National Basic Research Program of China (973 Program)  No.2007CB814800, the NSFC No.10631050 and No.10801084.

\end{document}